\documentclass[12pt]{article}
\usepackage{amsmath,amsthm,amssymb}
\usepackage{mathtools}

\newtheorem{theorem}{Theorem}[section] 
\newtheorem{lemma}{Lemma}[section]
\newtheorem{cor}{Corollary}[section]
\newtheorem{remark}{Remark}[section]

\newtheorem{proposition}{Proposition}[section]

\numberwithin{equation}{section}

\newcommand{\D}{{\rm d}}
\newcommand{\dx}{\, \D x}

\newcommand{\rz}{\mathbb{R}}
\newcommand{\nz}{\mathbb{N}}

\newcommand{\eps}{\varepsilon}
\newcommand{\iom}{\int_{\Omega}}

\newcommand{\klauf}{\left(\begin{array}}
\newcommand{\klzu}{\end{array}\right)}

\newcommand{\okappa}{{\hat{\kappa}}}

\newcommand{\ovsg}{\overline{s}_\alpha}
\newcommand{\ovs}{\overline{s}}
\newcommand{\sob}{\gamma}

\title{On the global regularity for minimizers of variational integrals: splitting-type problems in 2D
and extensions to the general anisotropic setting}

\author{Michael Bildhauer \& Martin Fuchs}
\date{}

\newcommand{\reff}[1]{(\ref{#1})}

\usepackage{hyperref}
\hypersetup{
   pdfnewwindow=true,
   colorlinks=false,
   linkcolor=black,
   citecolor=green,
   filecolor=magenta,
   urlcolor=cyan,
}

\newcommand{\olabel}[1]{\label{#1}}

\begin{document}

\maketitle

\newcommand{\op}[1]{\operatorname{#1}}

\newcommand{\hypref}[2]{\hyperref[#2]{#1 \ref*{#2}}}
\newcommand{\hypreff}[1]{\hyperref[#1]{(\ref*{#1})}}
\parindent0ex

\begin{abstract}
We mainly discuss superquadratic minimization problems for splitting-type variational integrals on a bounded Lipschitz domain 
$\Omega \subset \rz^2$ and prove
higher integrability of the gradient up to the boundary by incorporating an appropriate
weight-function measuring the distance of the solution to the boundary data. As a corollary, the local H\"older coefficient
with respect to some improved H\"older continuity is quantified in terms of the function $\op{dist}(\cdot,\partial \Omega)$.

The results are extended to anisotropic problems without splitting structure under
natural growth and ellipticity conditions.

In both cases we argue with variants of Caccioppoli's inequality involving small weights.%
{}
\end{abstract}

\newcommand{\qmin}{{q_{\min}}}
\newcommand{\qmax}{{q_{\max}}}

\parindent0ex
\section{Introduction}\label{intro}

In the main part of our note we consider splitting-type energies as a particular class of variational integrals
\begin{equation}\olabel{intro 1}
J[u] := \iom f(\nabla u) \dx
\end{equation}
defined for functions $u$: $\rz^n \supset \Omega \to \rz$.  In the splitting case, the energy density $f$: $\rz^n \to \rz$ admits
an additive decomposition into different parts depending on the various first partial derivatives of the admissible functions
$u$, for example, we can consider the case
\begin{equation}\olabel{intro 2}
f(\nabla u) = \sum_{i=1}^n f_i(\partial_i u)
\end{equation}
with functions $f_i$: $\rz \to \rz$, $i=1$, \dots , $n$, of possibly different growth rates.\\

This also leads us to the more general category of variational problems with non-standard growth, where in the simplest case one studies
an energy functional $J$ defined in equation \reff{intro 1} with convex density $f$ being bounded from above and from below
by different powers of $|\nabla u|$.\\
 
 For an overview of the various aspects of variational problems with non-standard growth including aspects as existence 
 and regularity we refer the interested reader, e.g., to \cite{Gi:1987_1}, \cite{Ma:1989_1}, \cite{Ch:1992_1}, \cite{ELM:1999_1}, 
 \cite{FS:1993_1}, \cite{Bi:1818}, \cite{BF:2007_1}, \cite{BFZ:2007_1} and to the recent paper \cite{BM:2020_1} together with the references
 quoted therein.\\ 
 
 Regarding the question of interior regularity of $J$-minimizers in the non-standard setting, the above references provide
 rather complete answers, while the problem of boundary regularity for the solution of
 \begin{equation}\olabel{intro 3}
 J[u] \to \min \, ,\qquad u=u_0 \quad\mbox{on}\; \partial \Omega\, ,
\end{equation}
with sufficiently smooth boundary datum $u_0$ seems to be not so well investigated.\\ 

We mention the contribution
\cite{LS:2014_1}, where the higher integrability of the minimizer $u$ up to the boundary is established for densities $f$
being in some sense close to the splitting-class defined in \reff{intro 2}. A precise formulation of the assumptions
concerning the density $f$ is given in inequality (2.1) of the paper \cite{LS:2014_1}, and we also like to stress that a survey
of the known global regularity results is provided in the introduction  of this reference.\\

Even more recently the contribution of Koch \cite{Ko:2021_1} addresses the global higher integrability of the gradient
of solutions of variational problems with $(p,q)$-growth, which to our knowledge is the first result to improve global
integrabilty of the gradient allowing full anisotropy with different growth rates with respect to different partial derivatives.
The boundary data are supposed to belong to some fractional order spaces (see, e.g., \cite{AF:2003_1} or \cite{DD:2012_1})
and roughly speaking are handled like an additional $x$-dependence following ideas as outlined in, e.g., \cite{ELM:2004_1}
for the interior situation. Hence, the admissible energy densities correspond to the examples presented in \cite{ELM:2004_1}.
On the other hand the case of vectorial functions $u$: $\rz^n \supset \Omega \to \rz^m$, 
$n\geq 2$, $m \geq 1$ as well as the subquadratic situation are included.\\

In our note we first discuss the variational problem \reff{intro 3} on a bounded Lipschitz domain $\Omega \subset \rz^2$ with energy functional
$J$ defined in equation \reff{intro 1} under the following assumptions on the data: for $i=1$, $2$ let $f_i \in C^2(\rz)$
satisfy
\begin{equation}\olabel{intro 4}
a_i (1+t^2)^{\frac{q_i-2}{2}} \leq f_i''(t) \leq A_i (1+t^2)^{\frac{q_i -2}{2}}\, , \qquad t\in \rz\, ,
\end{equation} 
with constants $a_i$, $A_i >0$ and exponents
\begin{equation}\olabel{intro 5}
q_i > 2 \, .
\end{equation}
We define the energy density $f$: $\rz^2 \to \rz$ according to (recall the ``splitting condition'' \reff{intro 2})
\begin{equation}\olabel{intro 6}
f(Z) := f_1(z_1) + f_2(z_2)\, ,\qquad Z = (z_1,z_2) \in \rz^2\, .
\end{equation}
Note that \reff{intro 4} yields the anisotropic ellipticity condition
\begin{equation}\olabel{intro 7}
c_1 (1+|Z|^2)^{\frac{p-2}{2}}|Y|^2  \leq D^2f(Z)(Y,Y) \leq c_2 (1+|Z|^2)^{\frac{q-2}{2}}|Y|^2
\end{equation}
valid for $Z$, $Y \in \rz^2$ with positive constants $c_1$, $c_2$ and for the choice of exponents $p=2$
and $q:= \max\{q_1,q_2\}$. In fact, \reff{intro 7} is a direct consequence of \reff{intro 5} and the
formula
\begin{equation}\olabel{intro 8}
D^2f(Z)(Y,Y) = \sum_{i=1}^2 f_i''(z_i)|y_i|^2\, ,\qquad Z, \, Y \in \rz^2 \, .
\end{equation}
Finally, we fix a functiom
\begin{equation}\olabel{intro 9}
u_0\in W^{1,\infty}(\Omega)
\end{equation}
remarking that actually $u_0 \in W^{1,t}(\Omega)$ with exponent $t$ being sufficiently large (depending on $q_1$ and $q_2$) is needed
in our calculations. For a definition of the Sobolev spaces $W^{k,r}(\Omega)$ and their local variants we refer to \cite{AF:2003_1}.\\

The function $u_0$ acts as a prescribed boundary datum in the minimization problem
\begin{equation}\olabel{intro 10}
J[u] := \iom f(\nabla u) \dx \to \min \qquad\mbox{in}\; u_0+W^{1,1}_0(\Omega)\, .
\end{equation}  
Now we state our result in the splitting case by the way establishing global higher integrabilty without relating $q_1$ and $q_2$
and under quite weak assumptions on $u_0$ which means that we do not have to impose additional hypotheses concerning the regularity of
the trace of $u_0$ in the sense that ${u_0}_{|\partial \Omega}$ belongs to some fractional Sobolev space.
The price we have to pay is a weight function which has to be incorporated.
\begin{theorem}\olabel{main}
Let the assumptions \reff{intro 4}-\reff{intro 6} and \reff{intro 9} hold. Moreover, let
$u \in u_0 +W^{1,1}_0(\Omega)$ denote the unique solution of problem \reff{intro 10}. Then it holds
\begin{equation}\olabel{intro 11}
|u-u_0|^t \Big(\big|\partial_1 u\big|^{\frac{3}{2}q_1} + \big|\partial_2 u\big|^{\frac{3}{2}q_2}\Big) \in L^1(\Omega) 
\end{equation}
for any choice of the exponent $t$ such that
\begin{equation}\olabel{intro 12}
t > T(q_1,q_2)
\end{equation}
with a finite lower bound $T(q_1,q_2) \geq 1$ depending on the growth rates $q_1$, $q_2$.
\end{theorem}
\begin{remark}\olabel{intro rem 1}
An admissible choice for the quantity $T(q_1,q_2)$ from \reff{intro 12} is 
\begin{equation}\olabel{intro 13}
T(q_1,q_2) = 6 \max\Bigg\{1, \frac{1}{2} \frac{\qmax}{\qmin}\frac{\qmin -1}{\qmin -2}, \frac{\qmin}{\qmin -2} - \frac{1}{6}\Bigg\}
\frac{\qmin}{\qmin -2} \, .
\end{equation}
Here and in the following we let 
\[
\qmin := \min\{q_1,q_2\}\, ,  \qquad \qmax := \max\{q_1,q_2\}\, .
\]
In the case $q_i >5$, $i=1$, $2$, we can replace \reff{intro 13} by the expression
\begin{equation}\olabel{intro 14}
T(q_1,q_2) = 6 \qmax \max\Bigg\{ \frac{1}{2\qmin -4},\frac{1}{\qmax+2}\Bigg\}
\end{equation}
which means that in the case of higher initial integrability we can decrease the exponent $t$ in \reff{intro 11}
leading to a stronger result. We also note that for the choice $q_i > 5$, $i=1$, $2$, the condition \reff{intro 9} can be replaced by
the requirement $u_0 \in W^{1,\qmin}(\Omega)$ together with $\partial_i u_0 \in L^{q_i}(\Omega)$, $i=1$, $2$.
\end{remark}

\begin{remark}\olabel{intro rem 2}
The existence and the uniqueness of a solution $u$ to problem \reff{intro 10} easily follow from the growth properties
of $f$ and its strict convexity combined with our assumption \reff{intro 9} on the boundary datum $u_0$.
Obviously the minimizer satisfies $\partial_i u\in L^{q_i}(\Omega)$, $i=1$, $2$, and from Theorem 1.1 ii) in \cite{BFZ:2007_1} we
deduce interior H\"older continuity of the first partial derivatives. This in turn implies that $u$ is in the
space $W^{2,2}_{\op{loc}}(\Omega)$.
\end{remark}

\begin{remark}\olabel{intro rem 3}
Of course we expect that our theorem can be extended to the case $n>2$, i.e.~to densities
\[
f(Z) = \sum_{i=1}^n f_i(z_i)\, \qquad Z \in \rz^n \, ,
\]
with functions $f_i$ satisfying \reff{intro 4} and \reff{intro 5} being replaced by $q_i > n$, $i=1$, \dots , $n$.
We restrict our considerations in the splitting case to the 2D-situation for two reasons: first, a generalization to the case $n \geq 2$
would be accompanied by a tremendous technical effort without promising a deeper insight. Secondly, the interior regularity results of \cite{BFZ:2007_1}
as the starting point for our arguments would have to be generalized to the full splitting case in $n$ dimensions which goes beyond the aim
of this note. 
\end{remark}

\begin{remark}\olabel{intro rem 4}
Inequality \reff{intro 11} can be seen as weighted global higher integrability result for $\nabla u$.
In fact, on account of \reff{intro 5}, $u$ is H\"older continuous up to the boundary, which means
$|u(x) - u_0(x)| \approx \op{dist}^\alpha (x,\partial \Omega)$ for points $x$ near $\partial \Omega$.
Thus, in this vague sense, Theorem \ref{main} yields $|\partial_i u|^{t_i} \op{dist}^\beta (\cdot , \partial \Omega) \in L^1(\Omega)$
for a certain range of exponents $\beta >0$ and $t_i > q_i$.\\

Another interpretation will be given in Corollary \ref{ein cor 1} on the behaviour of a suitable H\"older coefficient.
\end{remark}

In the second part of our paper we drop the splitting condition \reff{intro 6} and look at densities $f$ with anisotropic ($p$,$q$)-growth
in the sense of \reff{intro 7}
and for a moment we still restrict our considerations to the two-dimensional case. Then the arguments of the first part can be carried over
provided that $p$ and $q$ are not too far apart.
\begin{theorem}\label{theo nosplit}
Assume that $f$: $\rz^2 \to \rz$ is of class $C^2(\rz^2)$ satisfying \reff{intro 7} with exponents $2 < p \leq q$ such that
\begin{equation}\olabel{nosplit 1}
q < \min\Big\{p+\frac{1}{2},2p-2\Big\} \, .
\end{equation}
Let \reff{intro 9} hold for the boundary datum $u_0$ and let $u \in u_0 + W^{1,1}_0(\Omega)$ denote the unique solution
of problem \reff{intro 10}. Then we have ($\Gamma := 1 + |\nabla u|^2$)
\[
|u-u_0|^{2\kappa} \Gamma^s \in L^1(\Omega)
\]
for any exponents $\kappa$, $s$ such that
\begin{equation}\olabel{nosplit 2}
s\in \Big(\frac{q}{2}+\frac{3}{4},p-\frac{1}{4}\Big)
\end{equation}
and 
\begin{equation}\olabel{nosplit 3}
\kappa > \max\Bigg\{\frac{p}{p-2},\frac{s}{s- \big(\frac{q}{2}+\frac{3}{4}\Big)}, \frac{2s(p-1)}{(p-2)^2}\Bigg\} \, .
\end{equation}
\end{theorem}
\begin{remark}\olabel{intro rem 5}
Since $f$ is strictly convex and since we have \reff{intro 9}, existence and uniqueness of a solution $u$ to problem 
\reff{intro 10} are immediate. Note that \reff{nosplit 1} implies (4.2) in Marcellini's paper \cite{Ma:1991_1}, hence 
$u \in W^{1,q}_{\op{loc}}(\Omega)$ on account of Theorem 4.1 in \cite{Ma:1991_1}. Quoting Corollary 2.2 from this reference
we deduce $u \in C^{1,\alpha}(\Omega)$ and $u \in W^{2,2}_{\op{loc}}(\Omega)$.
\end{remark}

In our final theorem we incorporate two new features: we modify the idea of the proofs of the previous results
by starting with a more subtle inequality being applied to terms on the left-hand sides
of our foregoing calculations by the way extending the range of admissible exponents $p$ and $q$. 
Moreover, we include the general case of $n$-dimensions, $n \geq 2$. This can be done without technical efforts
comparable to the ones which are needed to generalize the 2D splitting case.\\

Now the precise formulation of our last result is as follows:
\begin{theorem}\label{aniso theo}
Consider the variational problem \reff{intro 10} on a bounded Lipschitz domain $\Omega \subset \rz^n$ now just 
assuming ($p$,$q$)-growth of the $C^2$-density $f$: $\rz^n \to \rz$ in the sense of \reff{intro 7} without splitting structure.
Let the boundary datum $u_0$ satisfy \reff{intro 9}.
Suppose that
\begin{equation}\olabel{aniso main 1}
2 \leq n < p \leq q \qquad\mbox{and}\qquad q < p + \frac{2 (p-n)}{n}\, .
\end{equation}
By assumption \reff{aniso main 1} we may choose a real number $\kappa > 1$ sufficiently large such that
\begin{equation}\olabel{aniso main 2}
q < p + \frac{2(p-n)}{n} -  \frac{2p-1}{n \kappa} \, .
\end{equation}
Moreover, let
\begin{equation}\olabel{aniso main 3}
\ovs < p \frac{\kappa(n+2)-2}{2n\kappa} - \frac{\kappa -1}{2n\kappa} \, .
\end{equation}
Increasing $\kappa$, if necessary, we suppose in addition that
\begin{equation}\olabel{aniso main 4}
\ovs < (\kappa -1) \frac{p-n}{n} + \frac{p}{2}
\end{equation}
and that we have the inequalities \reff{complete 7} below. Then it holds
\[
\iom \Gamma^{\ovs} |u-u_0|^{2(\kappa -1)} \dx \leq c \, , \qquad \Gamma := 1+|\nabla u|^2\, ,
\]
for a finite constant $c$ depending on $n$, $p$, $q$, $\kappa$ and $\ovs$.
\end{theorem}

\begin{remark}\label{main rem 1}
Passing to the limit $\kappa \to \infty$, \reff{aniso main 3} becomes
\[
\ovs \leq \frac{p}{2} + \frac{2p-1}{2n}\, .
\] 
\end{remark}

\begin{remark}\label{main rem 2}
We note that Remark \ref{intro rem 5} extends to the setting studied in Theorem \ref{aniso theo}, since now condition (4.2) from
\cite{Ma:1991_1} is a consequence of the more restrictive inequality \reff{aniso main 1}.
\end{remark}

\vspace*{2ex}
Our paper is organized as follows: in Section \ref{proof} we present the proof of Theorem \ref{main} for a given Lipschitz boundary datum 
$u_0$ (recall \reff{intro 9}) and just supposing \reff{intro 5}. In the case $q_1$, $q_2 > 5$ the lower bound $T(q_1,q_2)$
occurring in \reff{intro 12} looks much nicer and \reff{intro 9} can be considerably weakened. We sketch some arguments 
in Section \ref{fuenf}. The two-dimensional non-splitting case, i.e.~the setting described in Theorem \ref{theo nosplit}, is discussed in Section 
\ref{nosplit proof} and in Section \ref{aniso} we discuss the higher dimensional situation of Theorem \ref{aniso theo}.\\

During the proofs of our results we make essential use of Caccioppoli-type inequalities. As a main new feature these inequalities
involve ``small weights'', for instance weight functions like $(1+|\partial_i u|)^\alpha$, $i=1$, $2$, with a certain range
of negative exponents $\alpha$. A rather general analysis of these tools is presented in Section \ref{cacc}.\\

Before we get into all these technical details we finish the introduction by adding some extra comments concerning  the interpretation
of our results. For simplicity we restrict ourselves to the setting described in Theorem \ref{main} and denote by $u$
the unique solution of problem \reff{intro 10} under the hypothesis of Theorem \ref{main}.\\

We start with the following observation:

\begin{proposition}\label{ein prop 1}
Fix $x_0 \in \Omega$. Then there exists a neighbourhood $U=U(x_0)\Subset \Omega$ of $x_0$ such that we may assume that
for all $x \in U$
\[
|\hat{u}|(x) = |u-u_0| (x) \geq \underline{d} := \frac{1}{2} d(x_0) \, , \qquad \hat{u} := u-u_0\, ,
\]
where
\[
d(x) = \op{dist}(x,\partial \Omega) \, ,
\]
denotes the distance function to the boundary $\partial \Omega$.
\end{proposition}
{\it Proof of Proposition \ref{ein prop 1}.}  The distance function $d$ is discussed, e.g., in \cite{GT:1998_1}, Appendix 14.6. In particular,
if $\Omega$ is a domain in $\rz^n$ with non-empty boundary, then by (14.91) of \cite{GT:1998_1} $d$ is uniformly Lipschitz continuous.\\

Let us suppose w.l.o.g.~that $\hat{u}(x_0) \geq 0$.  If $\hat{u}(x_0) < d(x_0)$, then we choose $\tilde{u}_0 := u_0-d$ as admissible Lipschitz boundary datum 
with $u_0 = \tilde{u}_0$ on $\partial \Omega$.
Hence both $u_0$ and $\tilde{u}_0$ produce the same minimizer $u$ of the Dirichlet problem \reff{intro 10} and by definition we have
\[
(u-\tilde{u}_0)(x_0) = (u - u_0 + d)(x_0) \geq d(x_0) \, .
\]
The considerations of our note completely remain unchanged if we replace $u_0$ by $\tilde{u}_0$. Hence we have w.l.o.g.~that $\hat{u}(x_0) \geq d(x_0)$.
Proposition \ref{ein prop 1} then follows from the continuity of $u$ and $u_0$. \qed\\

Let us fix $x_0\in \Omega$, and define $U$, $\hat{u}$ as in Proposition \ref{ein prop 1}.  Moreover we let $h$: $\Omega \to \rz$,
$h:= |\hat{u}|^{\frac{4}{3} \frac{\kappa}{\qmin}+1}$,  hence
\[
|\nabla h| \leq \Bigg[\frac{4}{3} \frac{\kappa}{\qmin} +1\Bigg] |\hat{u}|^{\frac{4}{3}\frac{\kappa}{\qmin}}|\nabla \hat{u}|\, .
\]
Theorem \ref{main} gives
\[
\iom |\nabla h|^{\frac{3}{2}\qmin} \dx \leq c\, ,\qquad{i.e.}\qquad h\in W^{1,\frac{3}{2}\qmin}(\Omega) \, .
\]
Letting $\zeta = \frac{4}{3} \frac{\kappa}{\qmin}+1$, the imbedding into H\"older spaces yields $|\hat{u}|^{\zeta} \in C^{0,\mu}(\overline{\Omega})$
with $\mu = 1- \frac{4}{3\qmin}$.\\

Now we consider the Taylor-expansion of the function $g$: $\rz^+ \to \rz$, $w \mapsto w^{\zeta}$ around a fixed $\tilde{w} >0$: 
for $w$ sufficiently close to $\tilde{w}$ we have
\[
w^{\zeta} = \tilde{w}^\zeta + \zeta \tilde{w}^{\zeta -1} (w-\tilde{w}) + O(|w-\tilde{w}|^2) \, .
\]

Suppose w.l.o.g.~that $\hat{u} > 0$ in $U$. Inserting $\tilde{w} = \hat{u}(x_0)$ and $w = \hat{u}(x)$, $x$ sufficiently close to
$x_0$, in the Taylor-expansion we obtain using Proposition \ref{ein prop 1}
\begin{eqnarray}\label{ein 1}
\frac{|\hat{u}^\zeta(x) - \hat{u}^\zeta(x_0)|}{\zeta  \underline{d}^{\zeta -1}}& \geq & 
\frac{|\hat{u}^\zeta(x) - \hat{u}^\zeta(x_0)|}{\zeta \hat{u}^{\zeta -1}(x_0)}\nonumber\\
& \geq & u(x) - u(x_0) + u_0(x_0) - u_0(x) + O(|\hat{u}(x) - \hat{u}(x_0)|^2)\, .
\end{eqnarray}
Since $|\hat{u}|^\zeta$ is of class $C^{0,\mu}(\overline{\Omega})$ and since according to \reff{intro 9} $u_0$ is Lipschitz, 
by \reff{ein 1} we find a constant $c >0$ 
not depending on $x_0$ such that
\begin{equation}\label{ein 2}
u(x)- u(x_0) \leq c \frac{1}{\zeta \underline{d}^{\zeta -1}} |x-x_0|^\mu + O(|x-x_0|^{2\mu})\, .
\end{equation}
Changing the roles of $x$ and $x_0$ and following the arguments leading to \reff{ein 2}, we find a constant $c >0$ not depending on $x_0$
such that for all $x$ sufficiently close to $x_0$ we have
\begin{equation}
|u(x) - u(x_0)| \leq c \frac{1}{\zeta \underline{d}^{\zeta -1}}|x-x_0|^\mu \, .
\end{equation}

Summarizing the results we obtain the following corollary to Theorem \ref{main}.

\begin{cor}\label{ein cor 1}
Under the assumptions of Theorem \ref{main} we let 
\[
\zeta : =  \frac{4}{3} \frac{\kappa}{\qmin}+1 \, , \qquad \mu :=  1- \frac{4}{3\qmin} \, .
\]
If $x_0 \in \Omega$ is fixed, then for every sufficiently small neighbourhood $U = U(x_0)\Subset \Omega$ of $x_0$ 
the H\"older coefficient of $u$ (see \cite{GT:1998_1} for notation) satisfies
\[
[u]_{\mu;x_0} := \sup_{U} \frac{|u(x)-u(x_0)|}{|x-x_0|^\mu} \leq c \frac{1}{\zeta d^{\zeta-1}(x_0)} \, ,
\]
where the constant $c$ is not depending on $x_0$.
\end{cor}

\section{Proof of Theorem \ref{main}}\label{proof}


We proceed in several steps assuming from now on the validity of the hypotheses \reff{intro 4}, \reff{intro 5},
\reff{intro 6} and \reff{intro 9}. Let $u$ denote the unique solution of problem \reff{intro 10} and recall the
(interior) regularity properties of $u$ stated in \mbox{Remark \ref{intro rem 2}}. We start 
with the following elementary observation;

\begin{proposition}\olabel{proof prop 1}
Consider numbers $\kappa$, $s_1$, $s_2 \geq 1$. Then there is a finite constant $c$ such that for $i=1$, $2$ and
any $\eta \in C^1_0(\Omega)$
\begin{eqnarray}\olabel{proof 1}
\iom \Gamma_i^{s_i} |\eta|^{2\kappa} \dx &\leq &
c \Bigg\{ \Bigg[\iom \big| \nabla \Gamma_i^{\frac{s_i}{2}}\big|\,  |\eta|^\kappa \dx\Bigg]^2
+ \Bigg[\iom \Gamma_i^{\frac{s_i}{2}} |\nabla \eta| |\eta|^{\kappa-1} \dx\Bigg]^2\Bigg\}\nonumber\\
&=:& c \big\{ (I)_i +(II)_i\big\} \, ,\qquad \Gamma_i := 1+|\partial_i u|^2\, .
\end{eqnarray}
\end{proposition}

\emph{Proof of Proposition \ref{proof prop 1}.} Let $v:= \Gamma_i^{s_i/2} |\eta|^{\kappa}$
and observe that due to $u \in W^{2,2}_{\op{loc}} \cap W^{1,\infty}_{\op{loc}}(\Omega)$ the function $v$ is in the space $W^{1,1}_{0}(\Omega)$,
hence the Sobolev-Poincar\'{e} inequality implies 
\[
\iom v^2 \dx \leq c \Bigg[ \iom |\nabla v|\dx\Bigg]^2\, ,
\]
and inequality \reff{proof 1} directly follows by observing that $\big|\nabla |\eta|\big| \leq |\nabla \eta|$. \qed\\

Next  we replace $\eta$ by a sequence $\eta_m$, $m\in \nz$, being defined through
\begin{equation}\olabel{proof 2}
\eta_m := \varphi_m (u-u_0)
\end{equation}
with $\varphi_m \in C^1_0(\Omega)$, $0 \leq \varphi_m \leq 1$, 
\begin{eqnarray}\olabel{proof 3}
\varphi_m = 1 \; \mbox{on}\; \Big\{x\in \Omega: \, \op{dist}(x,\partial \Omega) \geq \frac{1}{m}\Big\}\, ,&&
\varphi_m = 0\; \mbox{on}\;  \Big\{x\in \Omega: \, \op{dist}(x,\partial \Omega) \leq \frac{1}{2m}\Big\}\, ,\nonumber\\
|\nabla \varphi_m| \leq c\,m \, ,&&  c=c(\partial \Omega)\, .
\end{eqnarray}
Note that due to \reff{intro 5} and \reff{intro 9} inequality \reff{proof 1} extends to $\eta_m$, moreover, according to
\cite{GT:1998_1}, Theorem 7.17, we have
\begin{equation}\olabel{proof 4}
|u(x) - u_0(x)| \leq c \rho^{1- \frac{2}{\qmin}}\, , \qquad \qmin := \min\{q_1,q_2\}\, ,
\end{equation}
for points $x\in \Omega$ at distance $\leq \rho$ to $\partial \Omega$, provided $\rho$ is sufficiently small.\\

The quantity $(II)_i$ defined in \reff{proof 1} with respect to the choice $\eta = \eta_m$ behaves as follows:
\begin{proposition}\olabel{proof prop 2}
Let $\delta := \qmin-2$ and consider numbers $s_i$, $\kappa \geq 1$, $i=1$, $2$, such that
\begin{equation}\olabel{proof 5}
s_i > \frac{\delta}{2}\, ,\qquad \kappa \geq 
\hat{\kappa} := \frac{\qmin -1}{\qmin-2}\, .
\end{equation}
Then it holds for $(II)_i = (II_m)_i$, $i=1$, $2$, defined in \reff{proof 1} with respect to the function $\eta_m$ from \reff{proof 2}
\begin{equation}\olabel{proof 6}
(II)_i \leq   c \iom \Gamma_i^{s_i - \frac{\delta}{2}} |\eta_m|^{2 \kappa - 2 \hat{\kappa}} \dx
=: c (III)_i \, , \qquad i=1 ,\, 2\, , 
\end{equation}
$c$ denoting a finite constant being uniform in $m$.
\end{proposition}
\emph{Proof of Proposition \ref{proof prop 2}.} By H\"older's inequality we have
\begin{eqnarray*}
(II)_i &=& \Bigg[ \iom \Gamma_i^{\frac{s_i}{2} - \frac{\delta}{4}} |\eta_m|^{\kappa -\okappa}
|\eta_m |^{\okappa-1} |\nabla \eta_m| \Gamma_i^{\frac{\delta}{4}}\dx\Bigg]^2\\
&\leq & \Bigg[ \iom \Gamma_i^{s_i - \frac{\delta}{2}} |\eta_m|^{2\kappa - 2 \okappa}\dx\Bigg]
\Bigg[ \iom \Gamma_i^{\frac{\delta}{2}} |\eta_m|^{2\okappa-2} |\nabla \eta_m|^2 \dx \Bigg]
\end{eqnarray*}
and for the second integral on the r.h.s.~we observe (recall \reff{proof 2} and \reff{proof 3}, \reff{proof 4})
\begin{eqnarray}\olabel{proof 6a}
|\eta_m|^{2\okappa-2} |\nabla \eta_m|^2 & \leq &
c \Big[ |\nabla (u-u_0)|^2 + |\nabla \varphi_m |^2 |u-u_0|^{2\okappa} \Big]\nonumber\\
&\leq & c \Big[ |\nabla u|^2 + |\nabla u_0|^2 + m^2 \Big(\frac{1}{m}\Big)^{2\okappa (1-2/\qmin)}\Big] \, .
\end{eqnarray}

Since $\Omega$ is Lipschitz we observe that  $|\op{spt} \nabla \varphi_m| \leq c/m$  and obtain
\begin{equation}\olabel{proof 6b}
\int_{\op{spt} \nabla \varphi_m} \Gamma_i^{\frac{\delta}{2}} m^{2-2\okappa (1-\frac{2}{\qmin})} \dx
\leq \iom \Gamma_i^{\frac{\qmin}{2}} \dx + m^{[2-2\okappa (1-\frac{2}{\qmin})]\frac{\qmin}{2}} m^{-1}
\end{equation}

The choice \reff{proof 5} of $\okappa$ shows that the right-hand side of \reff{proof 6b} is bounded.
Letting $\Gamma := 1 + |\nabla u|^2$ we arrive at
\begin{eqnarray*}
(II)_i&\leq & c\,  (III)_i  \Bigg[\iom \Big[ \Gamma^{\frac{\delta}{2}}\Gamma + \Gamma^{\frac{\delta}{2}} |\nabla u_0|^2\Big]\dx +1 \Bigg]\\
& \leq & c \, (III)_i \Bigg[ \iom \Gamma^{\frac{\qmin}{2}}\dx + \iom |\nabla u_0|^{\qmin}\dx +1 \Bigg] \,  ,
\end{eqnarray*}
where in the final step we have used the definition of $\delta$. Now \reff{proof 6} clearly is a consequence of \reff{intro 9} (in fact,
$u_0 \in W^{1,\qmin}(\Omega)$ would be sufficient).  \qed \\

In the following we will discuss the quantity $(III)_i$ defined in \reff{proof 6} under the assumptions of
Proposition \ref{proof prop 2}.\\

We have by Young's inequality for any $\eps > 0$
\begin{eqnarray*}
(III)_i &=& \iom \Gamma_i^{s_i \frac{s_i - \delta/2}{s_i}} |\eta_m|^{2 \kappa - 2 \hat{\kappa}}\dx
= \iom \Big[ \Gamma_i^{s_i} |\eta_m|^{2\kappa}\Big]^{\frac{s_i - \delta/2}{s_i}} 
|\eta_m|^{2\kappa - 2\hat{\kappa} - 2 \kappa \frac{s_i - \delta/2}{s_i}}\dx\\
&\leq &  \eps \iom \Gamma_i^{s_i} |\eta_m|^{2\kappa} \dx + c(\eps) \iom |\eta_m|^\vartheta \dx
\end{eqnarray*}
with exponent
\[
\vartheta := \Big[\frac{s_i}{s_i - \delta/2}\Big]^* \Big(2\kappa - 2 \hat{\kappa} - 2 \kappa \frac{s_i - \delta/2}{s_i}\Big) \, ,
\]
$[\dots ]^*$ denoting the exponent conjugate to $\frac{s_i}{s_i - \delta/2}$. Note that $\vartheta \geq 0$ provided we additionally assume that
$\kappa \geq 1$ satisfies the inequality
\begin{equation}\olabel{proof 7}
\kappa  \geq \frac{2s_i}{\delta} \hat{\kappa} = 2 s_i  \frac{\qmin -1}{(\qmin-2)^2}\, .
\end{equation}

Inserting the above estimate for $(III)_i$ into \reff{proof 6}, we find
\[
(II)_i \leq \eps \iom \Gamma_i^{s_i} |\eta_m|^{2\kappa} \dx + c(\eps) \, ,
\]
and if we choose $\eps$ sufficiently small, we see that \reff{proof 1} yields the following result:

\begin{proposition}\olabel{proof prop 3}
Let $s_i$, $\kappa \geq 1$ satisfy \reff{proof 5} and \reff{proof 7}. Then there exists a constant $c$ being independent of $m$ such that
($i=1$, $2$)
\begin{equation}\olabel{proof 8}
\iom \Gamma_i^{s_i} |\eta_m|^{2\kappa} \dx \leq  c \Big[ (I)_i+1\Big] \,  ,\qquad
(I)_i := \Bigg[\iom \big| \nabla \Gamma_i^{\frac{s_i}{2}}\big| \, |\eta_m|^\kappa \dx\Bigg]^2 \, . 
\end{equation}
\end{proposition}

In order to prove Theorem \ref{main} it remains to discuss the quantity $(I) = (I)_i$ for $i=1$, $2$, i.e.~now second
derivatives of $u$ have to be handled in an appropriate way. It holds
\begin{eqnarray*}
(I)_1& \leq & c \Bigg[ \iom \Gamma_1^{\frac{s_1}{2}-1} |\nabla \Gamma_1| \, |\eta_m|^{\kappa} \dx\Bigg]^2 \nonumber\\
&\leq & c \Bigg[ \iom \Gamma_1^{\frac{s_1-1}{2}} |\partial_1 \partial_1 u|\, |\eta_m|^{\kappa} \dx \Bigg]^2 +
c \Bigg[\iom \Gamma_1^{\frac{s_1 -1}{2}} |\partial_1 \partial_2 u|\, |\eta_m|^{\kappa}\dx\Bigg]^2
 =:  T_1 + \tilde{T}_1\, ,
\end{eqnarray*}
and in the same manner
\begin{eqnarray*}
(I)_2 &\leq & c \Bigg[ \iom \Gamma_2^{\frac{s_2-1}{2}} |\partial_2 \partial_2 u|\, |\eta_m|^{\kappa} \dx \Bigg]^2 +
c \Bigg[\iom \Gamma_2^{\frac{s_2 -1}{2}} |\partial_1 \partial_2 u|\, |\eta_m|^{\kappa}\dx\Bigg]^2
 =:  T_2 + \tilde{T}_2\, .
\end{eqnarray*}
We have
\begin{eqnarray*}
T_i &=& c \Bigg[ \iom \Gamma_i^{\frac{q_i-2}{4}} |\partial_i \partial_i u|\, |\eta_m|^{\kappa} \Gamma_i^{\frac{\alpha_i}{2}}
\Gamma_i^{\frac{s_i -1}{2}- \frac{q_i-2}{4} - \frac{\alpha_i}{2}} \dx \Bigg]^2\\
&  \leq & c \Bigg[ \iom \Gamma_i^{\frac{q_i-2}{2}}|\partial_i \partial_i u|^2 |\eta_m|^{2\kappa} \Gamma_i^{\alpha_i}\dx\Bigg]
\Bigg[ \iom \Gamma_i^{s_i -1 - \frac{q_i -2}{2} - \alpha_i} \dx \Bigg]\, ,
\end{eqnarray*}
where the last estimate follows from H\"older's inequality and $\alpha_i$, $i=1$, $2$, denote real numbers such that for the moment
\begin{equation}\olabel{proof 9}
s_i \leq q_i + \alpha _i \, , \qquad i=1, \, 2\, .
\end{equation}
Note that the condition \reff{proof 9} guarantees the validity of 
\begin{equation}\olabel{proof 10}
\iom \Gamma_i^{s_i - 1 - \frac{q_i -2}{2} - \alpha_i} \dx \leq c \, , \qquad i=1,\, 2\, ,
\end{equation}
for a finite constant $c$. Recalling \reff{intro 4} we see that \reff{proof 10} yields the bound
\begin{equation}\olabel{proof 11}
T_i \leq c \iom f_i''(\partial_i u) |\partial_i \partial_i u|^2 |\eta_m|^{2\kappa} \Gamma_i^{\alpha_i} \dx
\end{equation}
again for $i=1$, $2$.\\

Let us look at the quantities $\tilde{T}_i$: we have by H\"older's inequality
\begin{eqnarray}\olabel{proof 12}
\tilde{T}_1 &=& c \Bigg[ \iom \Gamma_1^{\frac{q_1-2}{4}} |\partial_1 \partial_2 u| |\eta_m|^{\kappa} \Gamma_1^{\frac{s_1-1}{2}} \Gamma_2^{\frac{\alpha}{2}}
\Gamma_1^{-\frac{q_1-2}{4}} \Gamma_2^{- \frac{\alpha_2}{2}}\dx\Bigg]^2\nonumber\\
& \leq & c \Bigg[\iom \Gamma_1^{\frac{q_1-2}{2}} |\partial_1 \partial_2 u|^2 |\eta_m|^{2\kappa} \Gamma_2^{\alpha_2} \dx \Bigg]
\Bigg[ \iom \Gamma_1^{s_1 -1 - \frac{q_1-2}{2}} \Gamma_2^{-\alpha_2}\dx\Bigg]\nonumber\\
&\leq & c \Bigg[ \iom f_1''(\partial_1 u) |\partial_1 \partial_2 u|^2 |\eta_m|^{2\kappa} \Gamma_2^{\alpha_2}\dx \Bigg]
\Bigg[\iom \Gamma_1^{s_1 -1 - \frac{q_1 -2}{2}} \Gamma_2^{-\alpha_2} \dx\Bigg]\nonumber\\
&=:& c S_1' \cdot S_1'' \, .
\end{eqnarray}
In order to benefit from the inequality \reff{proof 12}, we replace \reff{proof 9} by the stronger bound
\begin{equation}\olabel{proof 13}
s_i \leq \frac{3}{4} q_i \, ,\qquad i=1, \, 2\, ,
\end{equation}
together with the requirement
\begin{equation}\olabel{proof 14}
\alpha_i \in (-1/2,0)\, ,\qquad i=1, \, 2\, .
\end{equation}
Here we note that \reff{proof 13} together with \reff{proof 14} yields \reff{proof 9}. We then obtain on account of
Young's inequality
\[
S_1'' \leq c \Bigg[ \iom \Gamma_1^{\frac{q_1}{2}} \dx + \iom \Gamma_2^\beta \dx \Bigg]
\]
with exponent
\[
\beta := - \alpha_2 \Big[\frac{q_1/2}{s_1 - q_1/2}\Big]^* = - \alpha_2 \frac{q_1/2}{q_1-s_1} \leq - 2 \alpha _2 \leq \frac{q_2}{2} \, ,
\]
where we used \reff{proof 13}, \reff{proof 14} and \reff{proof 14}, \reff{intro 5}, respectively,  for the last inequalities. Thus \reff{proof 12}
reduces to
\begin{equation}\olabel{proof 15}
\tilde{T}_1 \leq c \iom f_1''(\partial_1 u) |\partial_1 \partial_2 u|^2 |\eta_m|^{2\kappa} \Gamma_2^{\alpha_2}\dx \, .
\end{equation}
In the same spirit it follows
\begin{equation}\olabel{proof 16}
\tilde{T}_2 \leq c \iom f_2''(\partial_2 u) |\partial_1 \partial_2 u|^2 |\eta_m|^{2\kappa} \Gamma_1^{\alpha_1}\dx \, .
\end{equation}
Let us return to \reff{proof 8}: we have by \reff{proof 11}, \reff{proof 15} and \reff{proof 16} and by \reff{intro 8}
\begin{eqnarray*}
\lefteqn{\iom \Big[\Gamma_1^{s_1} + \Gamma_2^{s_2}\Big] |\eta_m|^{2\kappa} \dx
\leq c \big[1+T_1+\tilde{T}_2+\tilde{T}_1+T_2]}\\
&=& c \Bigg[ 1+ \iom \Big(f_1''(\partial_1 u)|\partial_1\partial_1 u|^2 + f_2''(\partial_2 u) |\partial_1\partial_2 u|^2\Big)
|\eta_m|^{2\kappa} \Gamma_1^{\alpha_1} \dx\\
&&+  \iom \Big(f_1''(\partial_1 u)|\partial_1\partial_2 u|^2 + f_2''(\partial_2 u) |\partial_2\partial_2 u|^2\Big)
|\eta_m|^{2\kappa} \Gamma_2^{\alpha_2} \dx \Bigg]\\
&=& c \Bigg[ 1 + \iom D^2f(\nabla u) \big(\partial_1 \nabla u,\partial_1\nabla u\big) |\eta_m|^{2\kappa} \Gamma_1^{\alpha_1}\dx\\
&& + \iom D^2f(\nabla u) \big(\partial_2\nabla u,\partial_2 \nabla u\big)|\eta_m|^{2\kappa} \Gamma_2^{\alpha_2}\dx \Bigg] \, .
\end{eqnarray*}
The remaining integrals are handled with the help of Proposition \ref{prop cacc 1} (compare also inequality (4.6) from \cite{BF:2020_3}) 
replacing $l$ by $\kappa$ and $\eta$ by $|\eta_m|$, respectively. We note that the proof of Proposition \ref{prop cacc 1}
obviously remains valid with these replacements. We emphasize that \reff{proof 14} is an essential assumption to
apply Proposition \ref{prop cacc 1}. We get
\begin{eqnarray}\olabel{proof 17}
\iom \Big[\Gamma_1^{s_1} + \Gamma_2^{s_2}\Big] |\eta_m|^{2\kappa}
& \leq & 
c \Bigg[ 1 + \iom D^2f(\nabla u) \big(\nabla |\eta_m|,\nabla |\eta_m|\big) \Gamma_1^{1+\alpha_1}|\eta_m|^{2\kappa -2}\dx\nonumber\\
&& + \iom D^2f(\nabla u) \big(\nabla |\eta_m|,\nabla |\eta_m|\big) \Gamma_2^{1+\alpha_2} |\eta_m|^{2\kappa -2}\dx\Bigg]\, .
\end{eqnarray} 
Finally we let $s_i = \frac{3}{4}q_i$, $i=1$, $2$, which is the optimal choice with respect to \reff{proof 13}. We note that with this choice \reff{proof 7} 
follows from
\begin{equation}\olabel{proof 7a}
\kappa \geq \frac{3}{2} \frac{\qmax(\qmin-1)}{(\qmin -2)^2}\,  ,
\end{equation}
and \reff{proof 7a} is valid if $\kappa$ is chosen according to assumption \reff{intro 12} from Theorem \ref{main} with
$T(q_1,q_2)$ defined in \reff{intro 13}.\\

From \reff{intro 4} and \reff{intro 8} we deduce
\begin{eqnarray}\olabel{proof 18}
\mbox{r.h.s.~of \reff{proof 17}}& \leq & c \Bigg[ 1+ 
\sum_{i=1}^2 \iom \Gamma_i^{\frac{q_i}{2}+\alpha_i} |\partial_i \eta_m|^2 |\eta_m|^{2\kappa-2} \dx  \nonumber\\
&& + \iom \Gamma_2^{\frac{q_2 -2}{2}} |\partial_2 \eta_m|^2 |\eta_m|^{2\kappa -2} \Gamma_1^{1+\alpha_1}\dx\nonumber\\
&& + \iom \Gamma_1^{\frac{q_1-2}{2}} |\partial_1 \eta_m|^2 |\eta_m|^{2\kappa -2} \Gamma_2^{1+\alpha_2} \dx \Bigg]\, .
\end{eqnarray}
We recall the definition \reff{proof 2} of $\eta_m$ and the gradient bound for $\varphi_m$ stated in \reff{proof 3}. This
yields for integrand of the first integral on the right-hand side of \reff{proof 18}:
\begin{eqnarray*}
\Gamma_i^{\frac{q_i}{2}+\alpha_i} |\partial_i \eta_m|^2 |\eta_m|^{2 \kappa -2}
&\leq & c \Bigg[ \Gamma_i^{\frac{q_i}{2}+\alpha_i} |\partial_i (u-u_0)|^2 |\eta_m|^{2\kappa -2}
+\Gamma_i^{\frac{q_i}{2}+\alpha_i} |\partial_i \varphi_m|^2 |u-u_0|^{2\kappa}\Bigg]\, .
\end{eqnarray*}
We quote inequality \reff{proof 4} and recall that $\nabla \varphi_m$ has support in the set 
$\{x\in \Omega:\, \op{dist}(x,\partial\Omega)\leq 1/m\}$ satisfying $|\op{spt}\nabla \varphi_m|\leq c/m$, hence 
\begin{eqnarray*}
\iom \Gamma_i^{\frac{q_i}{2}+\alpha_i} |\partial_i \varphi_m|^2 |u-u_0|^{2\kappa} \dx
&\leq &  \int_{\op{spt}\nabla \varphi_m} \Gamma_i^{\frac{q_i}{2}+\alpha_i}m^{2-2\kappa \frac{\qmin-2}{\qmin}} \dx\\
&\leq & c \iom \Gamma_i^{\frac{q_i}{2}} \dx + c \int_{\op{spt} \nabla \varphi_m} m^{[ 2-2\kappa \frac{\qmin-2}{\qmin}] \gamma^*_i} \dx \, ,
\end{eqnarray*}
where (recall \reff{proof 14})
\[
\gamma_i = \frac{q_i}{q_i+2\alpha_i}\, ,\qquad \gamma_i^* = \Bigg[\frac{q_i}{q_i+2\alpha_i}\Bigg]^* = - \frac{q_i}{2\alpha_i}\, .
\]
Thus the second integral on the right-hand side is bounded if, with $\alpha_i$ sufficiently close to $-1/2$, we have
\begin{equation}\olabel{proof 18a}
\kappa > \frac{2\qmax-1}{2 \qmax} \frac{\qmin}{\qmin-2} \, .
\end{equation}
Assuming \reff{proof 18a} we arrive at
\begin{equation}\olabel{proof 19}
\iom \Gamma_i^{\frac{q_i}{2}+\alpha_i} |\partial_i \eta_m|^2 |\eta_m|^{2\kappa -2} \dx
\leq  c \Bigg[1 + \iom \Gamma_i^{\frac{q_i}{2}+\alpha_i}|\partial_i u - \partial_i u_0|^2 |\eta_m|^{2\kappa -2}\dx\Bigg]
\end{equation}
for any $\alpha_i$ sufficiently close to $-1/2$. Assuming this we next let
\[
\beta_i  = \frac{3q_i}{2q_i+4(1+\alpha_i)} \qquad\mbox{with conjugate exponent}\qquad 
\beta_i^*  = \frac{3q_i}{q_i - 4 (1+\alpha_i)}
\]
and apply Young's inequality in an obvious way to get
\[
\iom \Gamma_i^{\frac{q_i}{2}+\alpha_i+1} |\eta_m|^{2\kappa -2}\dx \leq \eps \iom \Gamma_i^{s_i}|\eta_m|^{2\kappa}\dx
+ c(\eps) \iom |\eta_m|^{2\kappa - 2 \beta_i^*}\dx \, ,
\] 
where the first term can be absorbed in the left-hand side of \reff{proof 17}, while the second one bounded under the assumption
\begin{equation} \olabel{proof 19a}
\kappa > \frac{3\qmin}{\qmin -2}\, .
\end{equation}
Here we used the fact that for $q > 2$ the function $q/(q-2)$ is a decreasing function. Altogether we have shown that for exponents $\alpha_i$ close
to $-1/2$ the first integral on the right-hand side of \reff{proof 18} splits into two parts, where the first one can be absorbed in the left-hand side
of \reff{proof 18} and the second one stays bounded. During our calculations we evidently used \reff{intro 9}, however \reff{intro 9}
can be replaced by weaker integrability assumtions concerning $\partial_i u$. We leave the details to the reader. \\

Let us finally consider the ``mixed terms'' on the right-hand side of \reff{proof 18}. We first observe the inequality
\begin{eqnarray}\olabel{proof 20}
\lefteqn{\iom \Gamma_2^{\frac{q_2-2}{2}} |\partial_2 \eta_m|^2 |\eta_m|^{2 \kappa -2}\Gamma_1^{1+\alpha_1}\dx}\nonumber\\
&\leq & c \iom \Gamma_2^{\frac{q_2-2}{2}} (\partial_2(u-u_0))^2   |\eta_m|^{2\kappa-2} \Gamma_1^{1+\alpha_1}\dx\nonumber\\
&&+\iom \Gamma_2^{\frac{q_2-2}{2}} |\partial_2 \varphi_m|^2 |u-u_0|^2 |\eta_m|^{2\kappa-2} \Gamma_1^{1+\alpha_1}\dx\, .
\end{eqnarray}
Considering the limit case $\alpha_1 =-1/2$, the first integral on the right-hand side of \reff{proof 20} basically is of the form
\[
\iom \Gamma_2^{\frac{q_2}{2}} \Gamma_1^{\frac{1}{2}} |\eta_m|^{2\kappa -2} \dx 
\]
and this integral directly results from an application of Caccioppoli's inequality. If we like to show an integrability result for
$\Gamma_2^{t_2}$ with some power $t_2 > q_2/2$, then the idea is to apply Young's inequality choosing
\[
\beta = \frac{2t_2}{q_2}\qquad\mbox{with conjugate exponent}\qquad \beta^*= \frac{2t_2}{2t_2-q_2}
\]
leading to the quantities $\Gamma_2^{t_2}$ and $\Gamma_1^{\frac{t_2}{2t_2 - q_2}}$.
If $t_1 =t_1(q_1) > q_1/2$ denotes the desired integrability exponent for $\Gamma_1$, then this requires the bound
\[
\frac{t_2}{2t_2-q_2} \leq t_1(q_1) \qquad\mbox{for all}\; q_1 > 2
\]
and of course we need the same condition changing the roles of $t_1$ and $t_2$.
With the symmetric Ansatz $t_i = \theta q_i$ for some $\theta > 1/2$ we are immediately led to $t_i = \frac{3}{4} q_i$,
$i=1$, $2$, which again motivates our choice of $s_i$.\\

More precisely: discussing the first integral on the right-hand side of \reff{proof 20} we choose $\beta_1 = 3/2$ with conjugate exponent
$\beta_1^* =3$ and obtain
\begin{eqnarray*}
\iom \Gamma_2^{\frac{q_2}{2}} |\eta_m|^{2\kappa -2} \Gamma_1^{1+\alpha_1} \dx 
& =  & \iom \Gamma_2^{\frac{q_2}{2}} \Gamma_1^{1+\alpha_1}  |\eta_m|^{\frac{2\kappa}{\beta_1}} |\eta_m|^{\frac{2\kappa}{\beta_1^*} -2} \dx \\
& \leq &\eps \iom \Gamma_2^{s_2} |\eta_m|^{2\kappa} \dx + c(\eps) \iom \Gamma_1^{3 (1+\alpha_1)} |\eta_m|^{2\kappa -6}\dx \, .
\end{eqnarray*}
Here the first integral is absorbed in the left-hand side of \reff{proof 17} and 
since $\alpha_1$ is chosen sufficiently close to $-1/2$, the second integral is bounded provided that we suppose in
addition
\begin{equation}\olabel{proof 21}
q_1 > 3 \quad\mbox{and}\quad \kappa \geq 3\, .
\end{equation}
If $q_1 < 3$, using
\[
\beta_2 =\frac{q_1}{4 (1+\alpha_1)} \, , \qquad \beta_2^* = \frac{q_1}{q_1-4(1+\alpha_1)}\, ,
\]
we are led to ($\tilde{\eps} \ll \eps$)
\[
c(\eps) \iom \Gamma_1^{3(1+\alpha_1)} \eta_m^{2\kappa -6} \dx \leq \tilde{\eps}\iom \Gamma_1^{s_1} \eta^{2\kappa} dx
+ c(\tilde{\eps}) \iom \eta^{2\kappa - 6 \beta_2^*}\dx\, ,
\]
where the first integral is absorbed in the right-hand side of \reff{proof 17} and
if $\alpha_1$ sufficiently close to $-1/2$ we now suppose
\begin{equation}\olabel{proof 22}
\kappa > \frac{3q_1}{q_1-2} \, .
\end{equation}
Note that the condition \reff{proof 22} is a consequence of \reff{proof 19a}.\\
 
It remains to discuss the second integral on the right-hand side of \reff{proof 20}. Using Young's inequality with
\[
\beta_3 = \frac{3}{2} \frac{q_2}{q_2 - 2}\, , \qquad \beta_3^* = \frac{3q_2}{q_2+4}\, ,
\]
we obtain
\begin{eqnarray}
\lefteqn{\iom \Gamma_2^{\frac{q_2-2}{2}} |\partial_2 \varphi_m|^2 |u-u_0|^2 |\eta_m|^{2\kappa-2} \Gamma_1^{1+\alpha_1}\dx}\nonumber\\
&\leq & \eps \iom \Gamma_2^{s_2} |\eta_m|^{2\kappa} \dx\nonumber\\ 
&&+ c(\eps) \iom \Gamma_1^{(1+\alpha_1) \frac{3q_2}{q_2+4}}
| \partial_2 \varphi_m|^{2\beta_3^*} |u-u_0|^{3 \beta_3^*}|\eta_m|^{2\kappa-2\beta_3^*}\dx\, ,
\end{eqnarray}
where again the first integral is absorbed in the left-hand side of \reff{proof 17}. Considering the second integral we choose
\[
\beta_4 = \frac{q_1 (q_2+4)}{4 q_2 (1+ \alpha_1)} \qquad\mbox{with conjugate exponent}\qquad
\beta_4^* = \frac{q_1(q_2+4)}{q_2[q_1-4(1+\alpha_1)]+4q_1} \, .
\]
This gives ($\tilde{\eps} \ll \eps$)
\begin{eqnarray}\olabel{proof 23}
 \lefteqn{c(\eps) \iom \Gamma_1^{(1+\alpha_1) \frac{3q_2}{q_2+4}}
| \partial_2 \varphi_m|^{2\beta_3^*} |u-u_0|^{3 \beta_3^*}|\eta_m|^{2\kappa-2\beta_3^*}\dx}\nonumber\\
&\leq & \tilde{\eps} \iom \Gamma_1^{s_1} |\eta_m|^{2\kappa}\dx \nonumber\\
&&+ c(\tilde{\eps}) \int_{\op{spt}\nabla \varphi_m} |\partial_2 \varphi_m|^{2 \beta_3^*\beta_4^*} |u-u_0|^{2\beta_3^* \beta_4^*}
|\eta_m|^{2\kappa-2\beta_3^*\beta_4^*}\dx\, .
\end{eqnarray}
As usual the first integral on the right-hand side of \reff{proof 23} is absorbed in the left-hand side of \reff{proof 17} and we calculate
\begin{equation}\olabel{proof 24}
\beta_3^* \beta_4^* = \frac{3 q_1 q_2}{q_2[q_1 - 4(1+\alpha_1)]+4q_1}\, .
\end{equation}
We obtain recalling \reff{proof 3} and \reff{proof 4}
\begin{eqnarray}
\lefteqn{ \int_{\op{spt}\nabla \varphi_m} |\partial_2 \varphi_m|^{2 \beta_3^*\beta_4^*} |u-u_0|^{2\beta_3^* \beta_4^*}
|\eta_m|^{2\kappa-2\beta_3^*\beta_4^*}\dx}\nonumber\\
& \leq &  \int_{\op{spt}\nabla \varphi_m} |\partial_2 \varphi_m|^{2 \beta_3^*\beta_4^*} |u-u_0|^{2\beta_3^* \beta_4^*}
|u-u_0|^{2\kappa-2\beta_3^*\beta_4^*}\dx\nonumber\\
& \leq &  c m^{-1} m^{\frac{6 q_1 q_2}{q_2[q_1 - 4(1+\alpha_1)]+4q_1}} m^{- \frac{\qmin-2}{\qmin} 2 \kappa} \, . 
\end{eqnarray}
For $\alpha_1$ sufficiently close to $-1/2$ this leads to the requirement 
\begin{equation}\olabel{proof 25}
\kappa  > 
 \frac{\qmin}{\qmin-2}\, \frac{1}{2} \Bigg[\frac{6q_1q_2}{q_2(q_1-2)+4q_1} - 1 \Bigg] \, .
\end{equation}
With $q_1>2$ fixed and for $q_2 > 2$ we consider the function $g(q_2) = \frac{6q_1q_2}{q_2(q_1-2) + 4q_1}$. We have
\begin{eqnarray*}
g'(q_2) &=& \frac{6 q_1}{q_2(q_1-2) + 4q_1} - \frac{6q_1q_2}{(q_2(q_1-2)+4q_1)^2} (q_1-2)\\
&=& \frac{24 q_1^2}{(q_2(q_1-2)+4q_1)^2} > 0 \, ,
\end{eqnarray*}
hence $g$ is an increasing function and
\[
g(q_2) \leq \lim_{t \to \infty} \frac{6 q_1 t}{t(q_1-2)+4q_1} = \frac{6q_1}{q_1-2} \, .
\]
Thus,
\begin{equation}\olabel{proof 26}
\kappa \geq \frac{\qmin}{\qmin-2}\, \frac{1}{2} \Bigg[6 \frac{\qmin}{\qmin-2} - 1\Bigg]
\end{equation}
implies the validity of \reff{proof 25}.\\

Summarizing the conditions imposed on $\kappa$ during our calculations, i.e.~recalling the bounds
\reff{proof 5}, \reff{proof 7a}, \reff{proof 18a}, \reff{proof 19a}, \reff{proof 26} we are led to the lower bound
\begin{equation}\olabel{proof 27}
\kappa > 3 \max\Bigg\{1, \frac{1}{2} \frac{\qmax}{\qmin}\frac{\qmin -1}{\qmin -2}, \frac{\qmin}{\qmin -2} - \frac{1}{6}\Bigg\} \frac{\qmin}{\qmin -2}
\end{equation}
for the exponent $\kappa$. Assuming the validity of \reff{proof 27} and returning to \reff{proof 17} we now have shown that for $\alpha_i$ sufficiently close
to $-1/2$ the right-hand side of \reff{proof 17} can be splitted into terms which either can be absorbed in the left-hand side
of \reff{proof 17} or stay uniformly bounded, hence
\begin{equation}\olabel{proof 28}
\iom \Big[\Gamma_1^{\frac{3}{4}q_1} + \Gamma_2^{\frac{3}{4}q_2}\Big] |\eta_m|^{2\kappa}\dx \leq c
\end{equation}
for a finite constant $c$ independent of $m$. Passing to the limit $m \to \infty$ in \reff{proof 28} our claim \reff{intro 11} follows.
Obviously \reff{intro 12} is a consequence of \reff{proof 27} and the definition of $T(q_1,q_2)$ stated in \reff{intro 13}, which completes the proof of
Theorem \ref{main} for arbitrary exponents $q_1$, $q_2 >2$. \qed

\begin{remark}\label{proof rem 1}
Let us add some comments on the behaviour of $\kappa$, which means that we look at the lower bound for the exponent $\kappa$
given by the right-hand side of inequality \reff{proof 27}.
\begin{enumerate}
\item Since $u-u_0 \in W^{1,\qmin}(\Omega)$ and thereby $u-u_0 \in C^{0,\nu}(\overline{\Omega})$, $\nu := 1- 2/\qmin$, the H\"older exponent
enters \reff{proof 27}, which also corresponds to the natural effect that
\item $\kappa \to \infty$ as $\qmin \to 2$.
\item The ratio $\qmax/\qmin$ determines the growth of $\kappa$.
\item In the limit $\qmax = \qmin \to \infty$ condition \reff{proof 27} reduces to $\kappa > 3$.
\end{enumerate}
\end{remark}

\section{Comments on the case $\qmin  > 5$}\label{fuenf}

We now choose the sequence $\eta_m \in C^1_0(\Omega)$ according to
\begin{eqnarray}\olabel{fuenf 1}
\partial_i \eta_m \to \partial_i u \, ,&&\mbox{in}\; L^{q_i}(\Omega) \, ,\quad  i= 1 ,\, 2 \, , \nonumber\\
\eta_m \to u-u_0 &&\mbox{uniformly as}\; m\to \infty \, .
\end{eqnarray}
For some elementary properties of anisotropic Sobolev spaces including an appropriate version of this density result we refer, e.g., 
to \cite{Ra:1979_1}, \cite{Ra:1981_1}). We emphasize
that during the following calculations condition \reff{intro 9} can be replaced by the weaker requirement
$\partial_i u_0 \in L^{q_i}(\Omega)$, $i=1$, $2$.
Proposition \ref{proof prop 1} obviously holds for $\eta_m$, and with $\delta$ as in Proposition \ref{proof prop 2} we now
obtain \reff{proof 6} with the choice $\okappa =1$ observing that going through the proof of \reff{proof 6} the quantity
$(II)_i$ can be handled as follows: first we note
\[
(II)_i \leq \Bigg[\iom \Gamma^{s_i - \frac{\delta}{2}} |\eta_m|^{2\kappa -2}\dx \Bigg]\cdot 
\Bigg[ \iom |\nabla \eta_m|^2 \Gamma_i^{\frac{\delta}{2}}\dx\Bigg]\, ,
\]
and then we use
\[
\iom |\nabla \eta_m|^2 \Gamma_i^{\frac{\delta}{2}} \dx \leq c \Bigg[1+\iom |\nabla u|^\qmin \dx\Bigg] \, ,
\]
which is a consequence of \reff{fuenf 1}.\\

We continue with the discussion of $(III)_i$ for the choice $\okappa =1$, and observe that \reff{proof 7} has to be replaced by
\begin{equation}\olabel{fuenf 2}
\kappa \geq \frac{2}{\delta} s_i\, .
\end{equation}
Exactly as in Section \ref{proof} we obtain the inequalities \reff{proof 17}, \reff{proof 18} now being valid for any exponents $s_i$, $\alpha_i$, $\kappa$
such that ($i=1$, $2$)
\begin{equation}\olabel{fuenf 3}
\alpha_i \in (-1/2,0)\, ,\qquad
s_i  \leq  \frac{3}{4} q_i\, ,\qquad
\kappa  \geq  \frac{2}{\delta} s_i \, .
\end{equation}
Still following the lines of Section \ref{proof} we let $s_i = \frac{3}{4} q_i$,  $i=1$, $2$ ,
and replace \reff{fuenf 2} by
\begin{equation}\olabel{fuenf 4}
\kappa \geq \frac{3}{2 \delta} \qmax\, . 
\end{equation}
We let
\[
\rho_i := 2 \kappa -2 -2\kappa \Big[\frac{q_i}{2}+\alpha_i\Big] \frac{1}{s_i}
\]
and consider the terms on the right-hand side of \reff{proof 18}.  Young's inequality yields ($0 < \eps < 1$)
\begin{eqnarray}\olabel{fuenf 5}
\lefteqn{\sum_{i=1}^2 \iom \Gamma_i^{\frac{q_i}{2}+\alpha_i} |\partial_i \eta_m|^2 |\eta_m|^{2\kappa-2} \dx}\nonumber\\
 &= & \sum_{i=1}^2 \iom\Big[\Gamma_i^{s_i} |\eta_m|^{2\kappa}\Big]^{\frac{1}{s_i} \big[\frac{q_i}{2}+\alpha_i\big]}
|\partial_i \eta_m|^2 |\eta_m|^{\rho_i}\dx\nonumber \\
&\leq & \sum_{i=1}^2 \Bigg[\eps \iom \Gamma_i^{s_i} |\eta_m|^{2 \kappa}\dx + c(\eps) \iom |\eta_m|^{\rho_i \gamma_i^*}
|\partial_i \eta_m|^{2 \gamma_i^*} \dx\Bigg] 
\end{eqnarray}
with
\[
\gamma_i := \frac{s_i}{\frac{q_i}{2}+\alpha_i} = \frac{3 q_i}{2 q_i + 4 \alpha_i}\, ,\qquad \gamma_i^* =\frac{3q_i}{q_i-4\alpha_i}\, .
\]
While the $\eps$-part of \reff{fuenf 5} can be absorbed, the remaining integrals are bounded if for $i=1$, $2$
\[
\rho_i \geq 0 \qquad\mbox{and}\qquad 2 \gamma_i^* \leq q_i \, .
\]
Noting that the function $t \mapsto t/(t+2)$, $t \geq 0$, is increasing and by choosing $\alpha_i$, $i=1$, $2$, sufficiently close to $-1/2$,
we see that these requirements are consequences of the strict inequalities
\begin{equation}\olabel{fuenf 6}
\qmax > 4\, ,\qquad \kappa  > \frac{3 \qmax}{\qmax +2} \, .
\end{equation}
We thus have assuming \reff{fuenf 4} and \reff{fuenf 6}
\begin{eqnarray}\olabel{fuenf 7}
\iom \Big[\Gamma_1^{s_1} + \Gamma_2^{s_2}\Big] |\eta_m|^{2\kappa}\dx
&\leq &  c \Bigg[ 1+ 
\iom \Gamma_2^{\frac{q_2 -2}{2}} |\partial_2 \eta_m|^2 |\eta_m|^{2\kappa -2} \Gamma_1^{1+\alpha_1}\dx\nonumber\\
&&+ \iom \Gamma_1^{\frac{q_1-2}{2}} |\partial_1 \eta_m|^2 |\eta_m|^{2\kappa -2} \Gamma_2^{1+\alpha_2} \dx \Bigg]\, .
\end{eqnarray}
Let us have closer look at the first integral on the right-hand side of \reff{fuenf 7}. Youngs's inequality gives
\begin{eqnarray}\olabel{fuenf 8}
\lefteqn{\iom \Gamma_2^{\frac{q_2 -2}{2}} |\partial_2 \eta_m|^2 |\eta_m|^{2\kappa -2} \Gamma_1^{1+\alpha_1}\dx}\nonumber\\
&\leq & \iom |\partial_2 \eta_m|^{q_2}\dx + 
\iom \Gamma_2^{\frac{q_2}{2}}  |\eta_m|^{(2\kappa -2) \frac{q_2}{q_2 -2}}
\Gamma_1^{(1+\alpha_1)\frac{q_2}{q_2 -2}}\dx\nonumber\\
&\leq & c\Bigg[1+ \iom \Big[ \Gamma_2^{s_2}|\eta_m|^{2\kappa}\Big]^{\frac{2}{3}} |\eta_m|^{(2\kappa-2) \frac{q_2}{q_2-2} - \frac{4}{3}\kappa}
\Gamma_1^{(1+\alpha_1) \frac{q_2}{q_2 -2}} \dx \Bigg]\nonumber\\
&\leq & c + \eps \iom \Gamma_2^{s_2} |\eta_m|^{2\kappa} \dx \nonumber\\
&&+ c(\eps) \iom |\eta_m|^{3 \big[(2\kappa -2) \frac{q_2}{q_2 -2} - \frac{4}{3} \kappa\big]} \Gamma_1^{3(1+\alpha_1)\frac{q_2}{q_2 -2}} \dx \, .
\end{eqnarray}
Here the first integral is absorbed and the second is bounded if we have (again $\alpha_2$ being sufficiently close to $-1/2$)
\[
(2\kappa-2) \frac{q_2}{q_2-2} - \frac{4}{3}\kappa \geq 0 \qquad\mbox{and}\qquad \frac{3}{2}\frac{q_2}{q_2-2} < \frac{q_1}{2}\, . 
\]
The first condition follows from the requirement \reff{fuenf 6},
the second one holds if we assume in addition that $\qmin > 5$.\\

In the same way the last term on the right-hand side of \reff{fuenf 7} is handled and 
by combining \reff{fuenf 4}, \reff{fuenf 6} we have shown Theorem \ref{main}
together with the formula \reff{intro 14}. \qed\\

\section{Proof of Theorem \ref{theo nosplit}}\label{nosplit proof}
We proceed along the lines of Section \ref{proof} assuming that all the hypothese of Theorem \ref{theo nosplit}
are satisfied. In place of \reff{proof 1} we have ($\Gamma:= 1+|\nabla u|^2$)
\begin{eqnarray}\olabel{pr 1}
\iom \Gamma^s |\eta|^{2\kappa} \dx &\leq & c \Bigg\{\Bigg[\iom |\nabla \Gamma^{\frac{s}{2}}| \, |\eta|^\kappa \dx\Bigg]^2
+\Bigg[ \iom \Gamma^{\frac{s}{2}} |\nabla \eta|\,|\eta|^{\kappa-1} \dx\Bigg]^2\Bigg\}\nonumber\\
& =: & c\{(I)+(II)\}\, ,
\end{eqnarray}
and \reff{pr 1} holds for any $s$, $\kappa \geq 1$ and $\eta \in C^1_0(\Omega)$. We define the sequence $\eta_m$ as in \reff{proof 2}
observing that \reff{proof 4} holds with $\qmin$ replaced by $p$. Letting $\delta := p-2$, defining
$\hat{\kappa} := \frac{p-1}{p-2}$ we obtain for
\begin{equation}\olabel{pr 2}
s > \frac{\delta}{2}\, , \quad \kappa \geq \hat{\kappa}
\end{equation}
applying obvious modifications in the proof of \reff{proof 6}
\begin{equation}\olabel{pr 3}
(II) \leq c \iom \Gamma^{s-\frac{\delta}{2}}|\eta_m|^{2\kappa - 2 \hat{\kappa}}\dx =: c (III) \, .
\end{equation}
Here $(II)$ just denotes the term occurring on the right-hand side of \reff{pr 1} with $\eta$ being replaced by $\eta_m$.
The quantity $(III)$ is discussed analogously to the term $(III)_i$: in accordance with the calculations presented after the proof of 
Proposition \ref{proof prop 2} we get for any $\eps >0$
\begin{eqnarray}\olabel{pr 4}
(III)&\leq & \eps \iom \Gamma^s |\eta_m|^{2\kappa}\dx +c(\eps) \iom |\eta_m|^\vartheta \dx\, ,\nonumber\\
\vartheta &:=& \Big[\frac{s}{s-\frac{\delta}{2}}\Big]^* \Big(2\kappa - 2 \hat{\kappa} - 2 \kappa \frac{s-\frac{\delta}{2}}{s}\Big)\, ,
\end{eqnarray}
wtith exponent $\vartheta \geq 0$, if we assume that
\begin{equation}\olabel{pr 5}
\kappa \geq 2 \frac{s}{\delta}\hat{\kappa}
\end{equation}
is satisfied. Choosing $\eps$ sufficiently small, inserting \reff{pr 4} in \reff{pr 3} and returning to \reff{pr 1} it is shown
(compare \reff{proof 9})
\begin{equation}\olabel{pr 6}
\iom \Gamma^s |\eta_m|^{2\kappa}\dx \leq c\big\{(I)+1\big\}\, \quad (I) := \Bigg[\iom |\nabla \Gamma^{\frac{s}{2}}| \,|\eta_m|^\kappa\dx\Bigg]^2\, .
\end{equation}
Adjusting the calculations presented after inequality \reff{proof 9} to the situation at hand we consider numbers $s$, $\alpha$
satisfying
\begin{equation}\olabel{pr 7}
s \leq p + \alpha
\end{equation}
and find recalling \reff{intro 7} and using \reff{pr 7}
\begin{eqnarray*}
(I) &\leq & c \Bigg[\iom \Gamma^{\frac{s-1}{2}}|\nabla^2 u|  |\eta_m|^\kappa \dx\Bigg]^2\\
&\leq & c \Bigg[ \iom \Gamma^{\frac{p-2}{2}} |\nabla^2 u|^2 |\eta_m|^{2\kappa} \Gamma^\alpha \dx\Bigg] \cdot
\Bigg[\iom \Gamma^{s-1-\frac{p-2}{2}-\alpha}\Bigg]\dx\\
&\leq & c \iom D^2f(\nabla u) \big(\partial_i \nabla u,\partial_i \nabla u\big) |\eta_m|^{2\kappa} \Gamma^\alpha \dx\, ,
\end{eqnarray*}
where from now on the sum is taken with respect to the index $i$. Let us assume in addition to \reff{pr 7} that we have
\begin{equation}\olabel{pr 8}
\alpha > - \frac{1}{4}
\end{equation}
as lower bound for the parameter $\alpha$. Quoting \reff{cacc 2} from Proposition \ref{prop cacc 2} (with $\eta$ being replaced by
$|\eta_m|^\kappa$) and returning to \reff{pr 6} we find:
\begin{equation}\olabel{pr 9}
\iom \Gamma^s |\eta_m|^{2\kappa} \dx \leq 
c \iom D^2f(\nabla u) \big(\nabla |\eta_m|,\nabla |\eta_m|\big) \Gamma^{1+\alpha} |\eta_m|^{2\kappa -2}\dx\, .
\end{equation}
On the right-hand side of \reff{pr 9} we use the second inequality from the ellipticity condition \reff{intro 7} to obtain
(in analogy to the inequalities \reff{proof 17}, \reff{proof 18}) 
\begin{equation}\olabel{pr 10}
\iom \Gamma^s |\eta_m|^{2\kappa} \dx \leq c \iom \Gamma^{\frac{q}{2}+\alpha}|\nabla \eta_m|^2 |\eta_m|^{2\kappa -2} \dx \, .
\end{equation}
The right-hand side of \reff{pr 10} is discussed following the arguments presented after \reff{proof 18}: we first have (recall \reff{proof 4})
\begin{eqnarray*}
\iom \Gamma^{\frac{q}{2}+\alpha} |\nabla \varphi_m|^2 |u-u_0|^{2\kappa}\dx &\leq &
\int_{\op{spt}\nabla \varphi_m} \Gamma^{\frac{q}{2}+\alpha} m^{2-2\kappa \frac{p-2}{p}} \dx\\
&\leq & c \iom \Gamma^{\frac{p}{2}}\dx + \int_{\op{spt}\nabla \varphi_m} m^{\gamma^*[2-2\kappa\frac{p-2}{p}]}\dx\, ,
\end{eqnarray*}
where we have defined
\[
\gamma:= \frac{p/2}{\frac{q}{2} +\alpha}\, , \qquad \gamma^* := \frac{\gamma}{\gamma-1}\, .
\]
Obviously this requires the bound $\frac{p}{2} \geq \frac{q}{2}+\alpha$, i.e.~the condition
\begin{equation}\olabel{pr 11}
q < p + \frac{1}{2}
\end{equation}
in combination with \reff{pr 8} yields the maximal range of anisotropy , since then we can choose $\alpha$ sufficiently close to 
$-1/4$ to guarantee $p > q+2\alpha$. Moreover, we assume the validity of $-1 + \gamma^*[2-2\kappa (p-2)/p]\leq 0$.
This is true again for $\alpha$ chosen sufficiently close to $-1/4$, provided that
\begin{equation}\olabel{pr 12}
\kappa > \frac{p+q-\frac{1}{2}}{2(p-2)}\, .
\end{equation}
Recalling \reff{pr 11} we see that \reff{pr 12} is a consequence of the stronger bound
\begin{equation}\olabel{pr 13}
\kappa > \frac{p}{p-2} \, .
\end{equation}
In accordance with \reff{proof 19} we therefore arrive at
\begin{equation}\olabel{pr 14}
\mbox{r.h.s.~of \reff{pr 10}} \leq c \Bigg[1+\iom \Gamma^{\frac{q}{2}+\alpha} |\nabla u - \nabla u_0|^2 |\eta_m|^{2\kappa -2}\dx \Bigg]\, .
\end{equation}
Neglecting the contribution resulting from $\nabla u_0$ and under the additional hypothesis
\begin{equation}\olabel{pr 15}
s > \frac{q}{2} + \frac{3}{4}\, ,
\end{equation}
which guarantees the validity of $s > \frac{q}{2}+1+\alpha$ for $\alpha > -1/4$, we have
\begin{eqnarray}\olabel{pr 16}
\iom \Gamma^{\frac{q}{2}+1+\alpha} |\eta_m|^{2\kappa -2}\dx &\leq & \eps 
\iom \Gamma^s |\eta_m|^{2\kappa} \dx + c(\eps) \iom |\eta_m|^{2\kappa - 2 \beta^*}\dx \, ,\\
\beta^* := \frac{\beta}{\beta-1}\, ,&& \beta := \frac{s}{1+\frac{q}{2}+\alpha}\, .
\end{eqnarray}
Let us add a comment concerning the conditions imposed on $s$: as remarked after \reff{pr 11}
and during the subsequent calculations the parameter $\alpha$ has to be adjusted in an appropriate way and might become very close
to the critical value $-1/4$. For this reason inequality \reff{pr 7} is replaced by the stronger one (recall \reff{nosplit 2})
\begin{equation}\olabel{pr 17}
s < p - \frac{1}{4}\, .
\end{equation}
In order to find numbers $s$ satisfying \reff{pr 15} and \reff{pr 17} we need the bound $q < 2p-2$ being more restrictive
for $p \in (2, 5/2)$ than inequality \reff{pr 11}.\\

Finally we note that the exponent $2\kappa - 2 \beta^*$ occurring in the second integral on the right-hand side of \reff{pr 16}
is non-negative provided $\kappa \geq s/(s-[1+\frac{q}{2}+\alpha])$, and the latter inequality holds for $\alpha$ sufficiently close to $-1/4$, if 
(compare \reff{nosplit 3})
\begin{equation}\olabel{pr 18}
\kappa > \frac{s}{s- \big[\frac{q}{2}+\frac{3}{4}\big]}\, .
\end{equation}
Assuming \reff{pr 18} our claim follows by inserting \reff{pr 16} into \reff{pr 14} and choosing $\eps$ sufficiently small. \qed

\section{Proof of Theorem \ref{aniso theo}}\label{aniso}

Let all the assumptions of Theorem \ref{aniso theo} hold. We define the number
\[
\sob := \frac{2n}{n+2}\, , \quad 1\leq \sob < 2\, , \qquad \mbox{with Sobolev conjugate}\qquad \frac{n\sob}{n-\sob} = 2 \, ,
\]
and note that
\[
\frac{n- \sob}{n} = \frac{n}{n+2} = \frac{\sob}{2}\, , \quad \frac{\sob}{2-\sob} = \frac{n}{2} \, .
\]
We also remark that for any fixed $\alpha$ and for any $\kappa >1$ the inequality
\begin{equation}\olabel{pre 3}
\frac{p}{2} < p \frac{\kappa(n+2)-2}{2n\kappa} +\alpha \frac{\kappa -1}{\kappa}\, .
\end{equation}
is equivalent to the requirement $\alpha > - p/n$.
\begin{lemma}\label{aniso lem 1}
Fix $\alpha > - 1/(2n)$, $\kappa > 1$ with \reff{pre 3} and choose a real number $\ovsg$ satisfying
\begin{equation}\olabel{aniso 1}
 \frac{p}{2} <  \ovsg < p \frac{\kappa(n+2)-2}{2n\kappa}
+ \alpha \frac{\kappa -1}{\kappa}\, .
\end{equation}

Then we have for any $\eps >0$ and for $\eta \in C^{1}_0(\Omega)$
\begin{equation}\olabel{aniso 2}
\iom \Gamma^{\ovsg} |\eta|^{2(\kappa -1)} \dx \leq 
\eps \iom \Gamma^{(\frac{q}{2}+\alpha)} |\nabla \eta|^{2} |\eta|^{2(\kappa -1)}\dx
+\eps  \iom \Gamma^{-\frac{p}{n} + s} |\nabla \eta|^{2} |\eta|^{2(\kappa -1)}\dx  + c
\end{equation}
for a constant being independent on $\eta$.
\end{lemma}

\emph{Proof of Lemma \ref{aniso lem 1}.}
We first observe that for
\begin{equation}\olabel{aniso 3}
s = \ovsg \frac{\kappa}{\kappa-1} - \frac{p}{2(\kappa -1)}
\end{equation}
we have after an application of Young's inequality
\begin{eqnarray}\olabel{aniso 4}
\iom \Gamma^{\ovsg} |\eta|^{2(\kappa-1)} \dx&=& \iom \Gamma^{\ovsg -\frac{p}{2\kappa}} |\eta|^{2(\kappa -1)} \Gamma^{\frac{p}{2\kappa}}\dx
\leq \tilde{\eps} \iom \Gamma^s |\eta|^{2 \kappa}\dx  + c(\tilde{\eps}) \iom \Gamma^{\frac{p}{2}}\dx 
\end{eqnarray}
for any $\tilde{\eps} > 0$. We then estimate using Sobolev's inequality
\begin{eqnarray}\olabel{aniso 5}
\iom \Gamma^s |\eta|^{2\kappa} \dx &=& \iom \big[\Gamma^{\frac{s}{2}} |\eta|^{\kappa}\big]^2\dx
\leq  c \Bigg[\iom \big| \nabla \big[ \Gamma^{\frac{s}{2}} |\eta|^\kappa \big]\big|^\sob \dx\Bigg]^{\frac{2}{\sob}}\nonumber\\
&\leq & c \Bigg[ \iom \big|\nabla \Gamma^{\frac{s}{2}}\big|^\sob |\eta|^{\kappa \sob}\dx\Bigg]^{\frac{2}{\sob}}
+ c(\kappa) \Bigg[\iom \Gamma^{\frac{s\sob}{2}} |\eta|^{(\kappa-1)\sob} |\nabla \eta|^{\sob} \dx\Bigg]^{\frac{2}{\sob}}\nonumber\\
&=& c T_1^{\frac{2}{\sob}} + c T_2^{\frac{2}{\sob}}\, .
\end{eqnarray}

Let us first consider $T_1$. H\"older's inequality gives
\begin{eqnarray}\olabel{aniso 6}
T_1^\frac{2}{\sob} &\leq & c \Bigg[ \iom |\nabla^2 u|^\sob \Gamma^{\sob \frac{s-1}{2}} |\eta|^{\kappa \sob} \dx\Bigg]^{\frac{2}{\sob}}\nonumber\\
&=& c \Bigg[ \iom |\nabla^2 u|^\sob \Gamma^{\sob \frac{p-2}{4}} \Gamma^{\sob\frac{\alpha}{2}}
 \Gamma^{\sob \frac{2-p}{4}} \Gamma^{-\sob\frac{\alpha}{2}}\Gamma^{\sob \frac{s-1}{2}} |\eta|^{\kappa \sob} \dx\Bigg]^{\frac{2}{\sob}}\nonumber\\
 &\leq & c \Bigg[ \iom |\nabla^2 u|^2 \Gamma^{\frac{p-2}{2}}\Gamma^{\alpha}|\eta|^{2\kappa}\dx\Bigg]
 \cdot \Bigg[ \iom \Gamma^{\frac{\sob}{2-\sob}( \frac{2-p}{2}-\alpha+s-1)}\dx\Bigg]^{\frac{2-\sob}{\sob}}
 \nonumber\\
&=& c T_{1,1} \cdot T_{1,2}^{\frac{2-\sob}{\sob}}\, .
\end{eqnarray}

For $T_{1,2}$ we observe recalling \reff{aniso 3}
\[
\frac{n}{2} \Big[ \frac{2-p}{2} -\alpha + s -1\Big] < \frac{p}{2}\Leftrightarrow \ovsg < 
p \frac{\kappa(n+2)-2}{2n \kappa} +\alpha\frac{\kappa -1}{\kappa}\, ,
\]
which is true by the choice \reff{aniso 1} of $\ovsg$, hence $T_{1,2}$ is uniformly bounded.\\

We handle $T_{1,1}$ with the help of Proposition \ref{prop cacc 2} (replacing $\eta$ by $|\eta|^\kappa$):
\begin{eqnarray}\olabel{aniso 7}
T_{1,1} & \leq & c \iom D^2f(\nabla u) \big(\nabla \partial_\sob u,\nabla \partial_\sob u\big) \Gamma^{\alpha} |\eta|^{2\kappa} \dx\nonumber\\
& \leq & c \iom D^2f(\nabla u)\big(\nabla |\eta|^\kappa , \nabla |\eta|^\kappa\big)|\nabla u|^2 \Gamma^{\alpha}\dx\nonumber\\
&\leq & c \iom \Gamma^{\frac{q-2}{2}} |\eta|^{2\kappa-2}|\nabla \eta|^2 \Gamma^{1+\alpha} \dx
= c \iom \Gamma^{\frac{q}{2} + \alpha} |\nabla \eta|^2 |\eta|^{2\kappa-2}\dx \, .
\end{eqnarray}

From \reff{aniso 5} - \reff{aniso 7} we conclude
\begin{equation}\olabel{aniso 8}
\iom \Gamma^s |\eta|^{2\kappa} \dx \leq c \iom \Gamma^{\frac{q}{2}-\gamma} |\nabla \eta|^2 |\eta|^{2(\kappa-1)}\dx 
+ c T_2^{\frac{2}{\sob}}
\end{equation}
with constants $c$ being independent of $\eta$ and
it remains to discuss $T_2$ in \reff{aniso 8}: we have for $\hat{\mu} =2/\sob$, $\mu = 2/(2-\sob)$,
\begin{eqnarray}\olabel{aniso 9}
T_2^{\frac{2}{\sob}} & = & \Bigg[\iom \Gamma^{\frac{s\sob}{2}}|\nabla \eta|^\sob |\eta|^{(\kappa-1)\sob}\dx\Bigg]^{\frac{2}{\sob}}
=  \Bigg[\iom \Gamma^{\frac{p}{2\mu}} \Gamma^{-\frac{p}{2\mu}}  
\Gamma^{\frac{s\sob}{2}}|\nabla \eta|^\sob |\eta|^{(\kappa-1)\sob}\dx\Bigg]^{\frac{2}{\sob}}\nonumber\\
&\leq & c \Bigg[\iom \Gamma^{\frac{p}{2}}\dx\Bigg]^{\frac{2}{n}} 
\cdot \iom \Gamma^{-\frac{p}{n}+s}|\nabla \eta|^{2}|\eta|^{2(\kappa-1)}\dx \, .
\end{eqnarray}

With \reff{aniso 4}, \reff{aniso 8} and \reff{aniso 9} the proof of Lemma \ref{aniso lem 1} is completed by choosing $\tilde{\eps}\ll \eps$. \qed\\

Now we come to the proof of the theorem: we choose $\eta_m = (u - u_0)\varphi_m$ with $\varphi_m$ defined after \reff{proof 2}.
Then Lemma \ref{aniso lem 1} yields
\begin{eqnarray}\olabel{complete 1}
\iom \Gamma^{\ovsg} |\eta_m|^{2(\kappa -1)} \dx
&\leq & \eps \iom \Gamma^{(\frac{q}{2}+\alpha)} \big[|\nabla u|+|\nabla u_0|\big]^{2} |\eta_m|^{2(\kappa -1)}\dx\nonumber\\
&&+\eps  \iom \Gamma^{-\frac{p}{n} + s} \big[|\nabla u| + |\nabla u_0|\big]^{2}  |\eta_m|^{2(\kappa -1)}\dx\nonumber\\
& & + \eps \iom \Gamma^{(\frac{q}{2}+\alpha)} |\nabla \varphi_m|^2 |u-u_0|^{2}|\eta_m|^{2(\kappa-1)}\dx\nonumber\\
 &&+\eps  \iom \Gamma^{-\frac{p}{n} + s} |\nabla \varphi_m|^2 |u-u_0|^2 |\eta_m|^{2(\kappa -1)}\dx + c \, .
\end{eqnarray}
Here the first two integrals on the right-hand side can be absorbed in the left-hand side provided that we have (recall that $u_0$ is Lipschitz)
\begin{equation}\olabel{complete 2}
\max\Big\{\frac{q}{2}+\alpha,  -\frac{p}{n} + s \Big\}<  \ovsg  -1 \, .
\end{equation}
Also on account of \reff{complete 2} we can handle the remaining integrals on the right-hand side of \reff{complete 1} with the
help of Young's inequality. We obtain
\begin{eqnarray}\olabel{complete 3}
\lefteqn{ \iom \Gamma^{(\frac{q}{2}+\alpha)} |\nabla \varphi_m|^2 |u-u_0|^{2}|\eta_m|^{2(\kappa-1)}\dx}\nonumber\\
 &\leq & \iom \Gamma^{\ovsg} |\eta_m|^{2(\kappa -1)} \dx
 + c \int_{\op{spt} \nabla \varphi_m} |\nabla \varphi_m|^{2\beta_1^*} |u-u_0|^{2\beta_1^*} |\eta_m|^{2(\kappa-1)} \dx \, .
\end{eqnarray}
as well as
\begin{eqnarray}\olabel{complete 4}
\lefteqn{ \iom \Gamma^{(-\frac{p}{n}+s)} |\nabla \varphi_m|^2 |u-u_0|^{2}\eta_m^{2(\kappa-1)}\dx}\nonumber\\
 &\leq & \iom \Gamma^{\ovsg} \eta_m^{2(\kappa -1)} \dx
 + c \int_{\op{spt} \nabla \varphi_m} |\nabla \varphi_m|^{2\beta_2^*} |u-u_0|^{2\beta_2^*} \eta_m^{2(\kappa-1)} \dx \, .
\end{eqnarray}
In \reff{complete 3} we choose ($\alpha$ sufficiently close to $-1/(2n)$) $\beta_1$ such that
\[
1 < \beta_1 < \ovsg \Big[\frac{q}{2}-\frac{1}{2n}\Big]^{-1} = 
\Bigg[p \frac{\kappa(n+2)-2}{\kappa} - \frac{\kappa -1}{\kappa}\Bigg] \frac{1}{qn -1}
\]
with conjugate exponent (recall $q-p < 2p/n$ on account of \reff{aniso main 1})
\begin{equation}\olabel{complete 5}
\beta_1^* > \frac{\kappa\big[p(n+2)-1\big]-2p+1}{\kappa\big[2p - (q-p)n\big]-2p+1}\, .
\end{equation}
In \reff{complete 4} we define $\beta_2$ according to
\[
1 < \beta_2 := \ovsg \Big[-\frac{p}{n}+s\Big]^{-1} = 
\Bigg[p \frac{\kappa(n+2)-2}{2\kappa} - \frac{\kappa -1}{2\kappa}\Bigg] \frac{1}{sn -p}
\]
with conjugate exponent
\begin{equation}\olabel{complete 6}
\beta_2^* > \frac{\kappa\big[pn +4p -1- 2sn\big]-2p+1}{\kappa\big[2p - (q-p)n\big]-2p+1} \, .
\end{equation}
Returning to \reff{complete 3} and \reff{complete 4}, respectively, the first integral again is absorbed in the left-hand side of \reff{complete 1}.
The remaining integrals stay bounded if we have for $i=1$, $2$ 
\[
m^{-1} m^{\beta_i^*\big[2 -2(1-n/p)\big]-2(\kappa -1)(1-n/p)} \leq c\, ,
\]
Thus we require the condition
\begin{equation}\olabel{complete 7}
\beta_i^* < (\kappa -1)\frac{p-n}{n} + \frac{p}{2n} \, , \qquad i=1,\, 2\, ,
\end{equation}
which on account of $p>n$ is satisfied for $\kappa$ sufficiently large.\\

It remains to arrange \reff{complete 2} together with the inequality on the right-hand side of \reff{aniso 1}. Here we first observe
\begin{eqnarray}\olabel{complete 8}
\frac{q}{2} + \alpha < p \frac{\kappa (n+2)-2}{2n\kappa} + \alpha \frac{\kappa -1}{\kappa} -1
&\Leftrightarrow & q < p \frac{\kappa (n+2)-2}{n \kappa} - 2 - 2 \frac{\alpha}{\kappa} \nonumber \\
&\Leftrightarrow & q-p < \frac{2}{n} (p-n) -  \frac{2p+2n\alpha}{n\kappa}
\end{eqnarray}
For $\alpha$ sufficiently close to $-1/(2n)$, \reff{complete 8} is a consequence of \reff{aniso main 2}.
We finally have to discuss 
\begin{eqnarray}\olabel{complete 9}
-\frac{p}{n} + s \leq \ovsg -1 &\Leftrightarrow & -\frac{p}{n} + \ovsg \frac{\kappa}{\kappa -1} - \frac{p}{2(\kappa -1)} < \ovsg - 1\nonumber\\
&\Leftrightarrow & \ovsg < (\kappa -1) \frac{p-n}{n} + \frac{p}{2}\, .
\end{eqnarray}
Assumption \reff{aniso main 4} implies \reff{complete 9}. Hence, \reff{complete 1}, \reff{complete 3} and \reff{complete 4} prove 
Theorem \ref{aniso} passing to the limit $m\to \infty$. \qed \\

\section{Appendix. Caccioppoli-type inequalities}\label{cacc}

We prove two Caccioppoli-type inequalities with small weights (i.e.~involving powers of $\Gamma = 1+ |\nabla u|^2$ or of
$\Gamma_i =1+|\partial_i u|^2$, $i=1$, \dots , $n$, with a certain range of negative exponents), where
the first one is the appropriate version in the splitting context.\\

There is no need to restrict the following considerations to the case $n=2$. Thus,
throughout this appendix, we suppose that $\Omega \subset \rz^n$ is a bounded Lipschitz domain and that 
$f$: $\rz^n \to \rz$ is of class $C^2$ satisfying $D^2f(Z)(Y,Y) > 0$ for all $Z$, $Y \in \rz^n$. \\

Moreover we suppose that $u \in W^{2,2}_{\op{loc}}(\Omega)\cap C^{1}(\Omega)$ solves the differentiated Euler equation 
\begin{equation}\olabel{cacc 1}
0 = \iom D^2 f(\nabla u) \big(\nabla \partial_i u, \nabla \psi\big) \dx \qquad\mbox{for all}\; \psi \in C^\infty_0(\Omega) 
\end{equation}
and for any $1 \leq i \leq n$ fixed.

\newcommand{\gam}[2]{\Gamma_{#1}^{#2}}
\begin{proposition}\label{prop cacc 1}
Fix $l\in \nz$ and suppose that $\eta \in C^\infty_0(\Omega)$, $0 \leq \eta \leq 1$.  
Then the inequality
\begin{eqnarray}\olabel{cacc 2}
\lefteqn{\iom D^2 f(\nabla u)\big(\nabla \partial_i u, \nabla \partial_i u\big) \eta^{2l} 
\gam{i}{\alpha} \dx}\nonumber \\ 
&& \leq c \iom D^2f(\nabla u) (\nabla \eta,\nabla \eta)\eta^{2l-2} \gam{i}{\alpha +1} \dx\, , \quad \Gamma_i := 1+|\partial_i u|^2\, ,
\end{eqnarray}
holds for any $\alpha > - 1/2$ and for any fixed $1\leq i \leq n$.
\end{proposition}

\emph{Proof of Proposition \ref{prop cacc 1}.} Suppose that $-1/2 < \alpha$ and fix $1 \leq i \leq n$ (no summation with respect to $i$).
Using approximation arguments we may insert 
\[
\psi := \eta^{2l} \partial_i u \gam{i}{\alpha}
\]
in the equation \reff{cacc 1} with the result
\begin{eqnarray}\olabel{cacc 3}
\iom D^2f(\nabla u) \big(\nabla \partial_i u , \nabla \partial_i u\big)
\eta^{2l} \gam{i}{\alpha}\dx
&=& - \iom D^2f(\nabla u)\big(\nabla \partial_i u, \nabla \gam{i}{\alpha}\big)
\partial_i u \eta^{2l}\dx\nonumber\\
&& - \iom D^2f(\nabla u)\big(\nabla \partial_i u, \nabla (\eta^{2l})\big)
\partial_i u \gam{i}{\alpha} \dx\nonumber\\
& =: & S_1+S_2 \, .
\end{eqnarray}
In \reff{cacc 3} we have 
\[
S_1 = - 2 \alpha \iom D^2f (\nabla u) \big(\nabla \partial_i u, \nabla \partial_i u\big) |\partial_i u|^2\gam{i}{\alpha-1}\eta^{2l}\dx 
\]
which gives $S_1 \leq 0$ if $\alpha \geq 0$. In this case we will just neglect $S_1$ in the following. In the case $-1/2 < \alpha < 0$ we estimate
\begin{eqnarray*}
|S_1| &=&  2 |\alpha| \iom D^2f (\nabla u) \big(\nabla \partial_i u, \nabla \partial_i u\big)
|\partial_i u|^2\gam{i}{\alpha-1}\eta^{2l}\dx\\
&\leq & 2  |\alpha| \iom D^2f (\nabla u) \big(\nabla \partial_i u, \nabla \partial_i u\big)\gam{i}{\alpha}\eta^{2l}\dx \,  .
\end{eqnarray*}
Since we have $2 |\alpha| < 1$ we may absorb $|S_1|$ in the left-hand side of \reff{cacc 3}, hence
\[
\iom D^2f (\nabla u) \big(\nabla \partial_i u, \nabla \partial_i u\big) \eta^{2l}\gam{i}{\alpha}\dx \leq c |S_2| \,  .
\]
 
For $0 < \eps$ sufficiently small we apply the Cauchy-Schwarz inequality to discuss $S_2$:
\begin{eqnarray*}
\lefteqn{\iom D^2f(\nabla u) (\nabla \partial_i u, \nabla \eta)\eta^{2l-1} 
\gam{i}{\alpha} \partial_iu \dx}\\
&\leq& \eps \iom D^2f(\nabla u) (\nabla \partial_i u ,\nabla \partial_i u) 
\eta^{2l} \gam{i}{\alpha}\dx\\ 
&&+ c(\eps) \iom D^2f(\nabla u) (\nabla \eta,\nabla \eta) 
\eta^{2l-2}\gam{i}{\alpha} |\partial_i u|^2 \dx\, .
\end{eqnarray*}
After absorbing the first term in the right-hand side of \reff{cacc 3} we have established our claim \reff{cacc 2}.  \qed\\

Instead of the quantities $\Gamma_i$ our second inequality (compare \cite{BF:2021_3} for the discussion in two dimensions)
involves the full derivative, i.e.~we incorporate certain negative powers of $\Gamma = 1+|\nabla u|^2$. As a consequence we do not obtain the range
$-1/2 < \alpha$ and have to replace this condition by the requirement $-1/(2n) < \alpha$.

\begin{proposition}\label{prop cacc 2}
Suppose that  $\eta\in C^{\infty}_{0}\big(\Omega\big)$ and fix some real number $\alpha$ such that $- 1/(2n) < \alpha$.
Then we have  (summation with respect to $i =1$, \dots , $n$)
\begin{eqnarray}\olabel{cacc 4}
\lefteqn{\big[1 + 2 \alpha n\big] \iom D^2 f (\nabla u) \big(\nabla \partial_i u, \nabla \partial_i u\big) 
\Gamma^{\alpha}\eta^2 \dx}\nonumber\\
&\leq & c \Bigg[ \int_{\op{spt}\nabla \eta} D^2 f (\nabla u) \big(\nabla \partial_i u, \nabla \partial_i u\big)
\Gamma^{\alpha}\eta^2 \dx\Bigg]^{\frac{1}{2}}
\nonumber\\
&& \qquad \cdot \Bigg[ \int_{\op{spt}\nabla \eta} D^2f(\nabla u)\big(\nabla \eta, \nabla \eta\big) \big|\nabla u\big|^2
\Gamma^{\alpha} \dx\Bigg]^{\frac{1}{2}}\, ,
\end{eqnarray}
where the constant $c$ is not depending on $\eta$.
In particular it holds
\begin{equation}\olabel{cacc 5}
\iom D^2 f (\nabla u) \big(\nabla \partial_i u, \nabla \partial_i u)
\Gamma^{\alpha} \eta^2\dx
\leq c  \int_{\op{spt}\nabla \eta} D^2f(\nabla u)\big(\nabla \eta, \nabla \eta\big) 
\big|\nabla u \big|^2\Gamma^{\alpha} \dx \, .
\end{equation}
\end{proposition} 

\emph{Proof.} For $i =1$, \dots , $n$ and any $\eta$ as above we deduce from \reff{cacc 1}
\begin{eqnarray}\olabel{cacc 7}
\lefteqn{\iom D^2f(\nabla u) \big(\nabla \partial_i u ,\nabla \partial_i u\big) 
\Gamma^{\alpha} \eta^2 \dx}\nonumber\\
&=& - \iom D^2f(\nabla u)\big(\nabla \partial_i u, \partial_i  u \nabla \Gamma^{\alpha}\big) \eta^2 \dx\nonumber\\
&& - 2 \iom D^2 f(\nabla u) \big(\nabla \partial_i u , \nabla \eta\big) \partial_i u \Gamma^{\alpha}\eta \dx 
=: I +II\, ,
\end{eqnarray}
where at this stage no summation with respect to the index $i$ is performed.
Following \cite{BF:2021_3} we denote the bilinear form $D^2f(\cdot,\cdot)$ by $\langle \cdot,\cdot\rangle$. The first observation
is the inequality
\begin{eqnarray}\olabel{cacc 8}
\sum_{i =1}^n \Big\langle \partial_i \nabla u, \partial_i \nabla u\Big\rangle \Gamma^{\alpha}
&\geq & \sum_{i =1}^n \Big\langle \partial_i \nabla u, \partial_i \nabla u\Big\rangle
\sum_{j=1}^n (\partial_j u)^2   \Gamma^{\alpha -1}\nonumber \\
&\geq & \sum_{i =1}^n \Big\langle \partial_i \nabla u, \partial_i \nabla u\Big\rangle
(\partial_i u)^2   \Gamma^{\alpha -1}\nonumber \\
& = & \sum_{i =1}^n \Big\langle \partial_{i} u  \partial_i \nabla u , 
\partial_i u   \partial_i \nabla u \Big\rangle \Gamma^{\alpha -1} \, .
\end{eqnarray}

For handling the integrand of the term $I$ from \reff{cacc 7} we use the identity
\begin{eqnarray}\label{cacc 9}
- \sum_{i=1}^n\Big\langle \partial_i \nabla u , \partial_i u \nabla  \Gamma^{\alpha}\Big\rangle
&=& - \alpha \sum_{i=1}^n \Big\langle \partial_i u  \partial_i \nabla u , 
\nabla \sum_{j=1}^n (\partial_j u )^2\Big\rangle \Gamma^{\alpha -1}\nonumber\\ 
 &=& - 2\alpha \sum_{i =1}^n \Big\langle \partial_i u  \partial_i \nabla u,  
 \sum_{j =1}^n \partial_j u \partial_j \nabla u\Big\rangle \Gamma^{\alpha -1}\nonumber\\
 &=& -2 \alpha \sum_{i =1}^n \Big\langle \partial_i u  \partial_i \nabla u,  \partial_i u \partial_i \nabla u\Big\rangle
 \Gamma^{\alpha -1}\nonumber\\
 && - 2 \alpha \sum_{i =1}^n \Big\langle \partial_i u  \partial_i \nabla u,  
 \sum_{j\not=i}\partial_j u \partial_j \nabla u\Big\rangle
 \Gamma^{\alpha -1} \, .
\end{eqnarray}

The last term on the right-hand side of \reff{cacc 9} is estimated as follows
\begin{eqnarray*}
\lefteqn{\sum_{i =1}^n \Big\langle \partial_i u  \partial_i \nabla u,  
 \sum_{j\not=i}\partial_j u \partial_j \nabla u\Big\rangle}\\
 &\leq & \sum_{i=1}^n  \Bigg[\sum_{j\not= i}
 \frac{1}{2} \Big[ \Big\langle \partial_i u  \partial_i \nabla u,  \partial_i u  \partial_i \nabla u\Big\rangle 
 +\Big\langle \partial_j u \partial_j \nabla u,  \partial_j u \partial_j \nabla u\Big\rangle\Big]\Bigg]\\
 &=& \frac{1}{2} \sum_{i =1}^n 
 \Bigg[ (n-1) \Big\langle \partial_i u  \partial_i \nabla u,  \partial_i u  \partial_i \nabla u\Big\rangle 
 + \sum_{j\not=i} \Big\langle \partial_j u \partial_j \nabla u,  \partial_j u \partial_j \nabla u\Big\rangle\Big]\Bigg]\\
 &=& \frac{1}{2} \sum_{i =1}^n 
 \Bigg[ (n-2) \Big\langle \partial_i u  \partial_i \nabla u,  \partial_i u  \partial_i \nabla u\Big\rangle 
 + \sum_{j=1}^n \Big\langle \partial_j u \partial_j \nabla u,  \partial_j u \partial_j \nabla u\Big\rangle\Big]\Bigg]\\
 &=& \frac{1}{2} \Bigg[ (n-2) \sum_{i =1}^n 
 \Big\langle \partial_i u  \partial_i \nabla u,  \partial_i u  \partial_i \nabla u\Big\rangle 
 + n  \sum_{j=1}^n \Big\langle \partial_j u \partial_j \nabla u,  \partial_j u \partial_j \nabla u\Big\rangle\Big]\Bigg]\\
 &=& (n-1) \sum_{i=1}^n  \Big\langle \partial_i u  \partial_i \nabla u,  \partial_i u  \partial_i \nabla u\Big\rangle \, .
\end{eqnarray*}

This, together with \reff{cacc 9} gives (recall $-2\alpha > 0$)
\begin{eqnarray}\label{cacc 10}
- \sum_{i=1}^n\Big\langle \partial_i \nabla u , \partial_i u \nabla  \Gamma^{\alpha}\Big\rangle &=&
 -2 \alpha \sum_{i =1}^n \Big\langle \partial_i u  \partial_i \nabla u,  \partial_i u \partial_i \nabla u\Big\rangle
 \Gamma^{\alpha -1}\nonumber\\
 && - 2 \alpha \sum_{i =1}^n \Big\langle \partial_i u  \partial_i \nabla u,  
 \sum_{j\not=i}\partial_j u \partial_j \nabla u\Big\rangle
 \Gamma^{\alpha -1} \nonumber\\
&\leq &   -2 \alpha \sum_{i =1}^n \Big\langle \partial_i u  \partial_i \nabla u,  \partial_i u \partial_i \nabla u\Big\rangle
 \Gamma^{\alpha -1}\nonumber\\
&&- 2 \alpha (n-1) \sum_{i =1}^n \Big\langle \partial_i u  \partial_i \nabla u, 
\partial_i u  \partial_i \nabla u\Big\rangle
 \Gamma^{\alpha -1} \nonumber\\
 &=&- 2 \alpha n \sum_{i =1}^n \Big\langle \partial_i u  \partial_i \nabla u, 
\partial_i u  \partial_i \nabla u\Big\rangle
 \Gamma^{\alpha -1} \, .
\end{eqnarray}

Combining \reff{cacc 10} and \reff{cacc 8} we get
\begin{eqnarray}\label{cacc 11}
\lefteqn{-\sum_{i =1}^n  \iom D^2f(\nabla u)\big( \partial_i \nabla u , 
\partial_i u  \nabla  \Gamma^{\alpha}\big)\eta^2 \dx} \nonumber\\
&&\leq -2 \alpha n \sum_{i =1}^n
\iom D^2f(\nabla u) \big(\partial_i \nabla u, \partial_i \nabla u\big)\Gamma^{\alpha} \eta^2 \dx\, .
\end{eqnarray}

Returning to \reff{cacc 7} and using \reff{cacc 11} we get (from now on summation with respect to $i$)
\begin{eqnarray}\label{cacc 12}
\lefteqn{\big[1+ 2 \alpha n \big]  \iom D^2 f(\nabla u) \big(\nabla \partial_i u , \nabla \partial_i u\big) 
\Gamma^{\alpha} \eta^2 \dx}\nonumber\\
&&  \leq -2 \int_{\op{spt}\nabla \eta} D^2 f(\nabla u) \big(\eta \nabla \partial_i u, \partial_i u   \nabla \eta\big) 
\Gamma^{\alpha} \dx \, .
\end{eqnarray}

This finishes the proof by applying the Cauchy-Schwarz inequality. \qed

\bibliography{global_int}
\bibliographystyle{unsrt}

\end{document}